\documentclass[12 pt]{article}
\usepackage[english]{babel}
\usepackage[active]{srcltx}
\usepackage{amsmath}
\usepackage{amsfonts,amssymb}
\usepackage[mathcal]{eucal}

\newcommand{\rav}{\triangleq}
\newcommand{\av}{v}
\newcommand{\bw}{w}
\newcommand{\avm}{\mathcal{V}^{\,\flat}}
\newcommand{\avp}{\mathcal{V}^{\,\sharp}}
\newcommand{\bwm}{\mathcal{W}^{\,\flat}}
\newcommand{\bwp}{\mathcal{W}^{\,\sharp}}
\newcommand{\cyp}{{\mathbb{V}}^{\,\!\sharp}}
\newcommand{\cym}{{\mathbb{V}}^{\,\!\flat}}
\newcommand{\avg}{\mathcal{V}^{\,\natural}}
\newcommand{\bwg}{\mathcal{W}^{\,\natural}}
\newcommand{\bavg}{\overline{\mathcal{V}}^{\, }}
\newcommand{\bbwg}{\overline{\mathcal{W}}^{\, }}

\newcommand{\rref}[1]{$(\ref{#1})$}
\newcommand{\pust}{\varnothing}
\newcommand{\mm}[1]{{\mathbb{#1}}}
\newcommand{\mA}{{\mathfrak{A}}}
\newcommand{\mB}{{\mathfrak{B}}}
\newcommand{\fr}[1]{{\mathfrak{#1}}}
\newcommand{\ct}[1]{{\mathcal{#1}}}
\newcommand{\ravref}[1]{\stackrel{(\ref{#1})}{=}}
\newcommand{\leqref}[1]{\stackrel{(\ref{#1})}{\leq}}
\newcommand{\geqref}[1]{\stackrel{(\ref{#1})}{\geq}}
\newcommand{\epsi}{\varepsilon}
\newcommand{\bo}{\hfill {$\Box$}}
\newcommand{\ldr}{\langle}
\newcommand{\rdr}{\rangle}
\newcommand{\proof}{{\it Proof }}
\newtheorem{definition}{Definition}
\newtheorem{proposition}{Proposition}
\newtheorem{remark}{Remark}
\newtheorem{theorem}{Theorem}
\newtheorem{lemma}{Lemma}

\newtheorem{corollary}{Corollary}
\renewcommand{\baselinestretch}{1.3}
\begin{document}

\title{On  Uniform Tauberian Theorems for Dynamic Games\thanks{
               Krasovskii Institute,, Yekaterinburg, Russia;\ \  Ural Federal University, Yekaterinburg, Russia
}
}

\author{Dmitry Khlopin\\
{\it khlopin@imm.uran.ru} 
}

\maketitle

\begin{abstract}
The paper is concerned with two-person dynamic zero-sum games in continuous setting.
    We investigate the limit of value functions of finite horizon games with long run average cost as the time horizon tends
  to infinity and the limit of value functions of $\lambda$-discounted games as the discount tends to zero.
      Under quite weak assumptions on the game, we prove the  Uniform Tauberian Theorem: existence  of a uniform limit for one of the value functions implies
      the uniform convergence of the other one to the same limit.
      We also prove the analogs of the One-sided Tauberian Theorem, i.e., the inequalities on asymptotics for the lower and upper games. Also, a variant of the theorem for discrete-time games is treated separately.
  Special attention is devoted to the case of differential games.  The key roles in the proof were played by Bellman's optimality principle and
    the closedness of strategies under concatenation.

\end{abstract}

 Hardy once proved (see, for example,  \cite{Hardy1914})
  that, for a bounded sequence of numbers $a_i$,
 $$\lim_{n\to\infty} \frac{1}{n}\sum_{i=1}^n a_i=
 \lim_{\lambda\downarrow 0} \lambda \sum_{i=1}^\infty (1-\lambda)^{i-1} a_i$$
 if there exists at least one of these limits.
  This result was named a Tauberian theorem in honor of the similar result obtained by Tauber for convergent series. Theorems of this kind, in particular, provide the means for obtaining good estimates for sums of series with the use of faster methods of summation. For a more detailed treatment of the history of those results, see, for example, \cite{Bingham}.
   There are also Tauberian theorems for functions, in particular (see, for example, \cite[Sect. 6.8]{Hardy1949}), for a bounded continuous function $g$, the limit of  long run averages  and  limit of discounted averages (Cesaro mean  and Abel mean, respectively)
  $$\frac{1}{T} \int_{0}^T g(t)\,dt,\qquad
    \lambda \int_{0}^\infty e^{-\lambda t}g(t)\,dt$$
    coincide if there exists at least one of these limits.

What if we optimize the Abel mean and/or Cesaro mean and then consider the limit of the optimal values corresponding to them? Such limit value (as the discount tends to zero) was first considered in
  \cite{Kh_Blackwell} for a stochastic formulation. As proved in
 \cite{MN1981}, for a stochastic two-person game with a finite number of states and actions, optimal long-time~averages and optimal discounted averages share the common limit. For more details on  the limit value for Abel mean and/or Cesaro mean in other stochastic formulations, see \cite{BGQ2013,Vigeral,Ziliotto}.

 In the deterministic case, the question of existence of limit values arose in the control theory, time and again; one may at the very least note
\cite{chl1991,Gaitsgori1986,grune,Lions_unpub}.
 In the ergodic case (more generally, in the nonexpansive-like case) such limits exist and, moreover, they are independent of  the initial state
\cite{AlvBardi,ArisawaLions,AG2010,AG2000,BardiGame}.
For the latest results on existence of the limit values (first of all, in the nonergodic case), refer to
\cite{QG2009,QG2013,QuinRenault}; see also the review in  \cite[Sect. 3.4]{BCQ2013}.
  For a bibliography on discrete statements, refer to \cite{Renault2011}.

For control problems, the equality of limit values in the case when at least one of those limits is constant was first proved in  \cite{arisconst}.
 The very general case of dynamic system  was considered in paper \cite{barton}.
 Namely, it was proved there that existence of a uniform (on the set of states) limit for the value of one of the means implies the uniform convergence for value of the other one to the same limit. That paper  also contains  a beautiful introduction to the history of the subject and a review of publications in the field.  For discrete time systems, the equality of limit values was proved earlier in  \cite{Lehrer}.

Until very recently,
there was only a rather small number of publications concerned with the study of the limits of optimal averages in differential games. These are, first of all,  \cite{AlvBarditrue,BardiGame,Cardal2010} and, in addition, \cite{Bettiol,Ghosh}.
 In differential games of the special kind, those limits may be connected with the  asymptotic value of zero-sum repeated games; a good discussion of this issue is presented in
 \cite{Sorin},  see \cite{SIAM2012} on Tauberian theorem for repeated games. Moreover,
 as noted in \cite{barton} for differential games,
 ``When
  the dynamic is controlled by two players with opposite goals,
   a Tauberian theorem is given in
  the ergodic case by Theorem 2.1
  in \cite{AlvBarditrue}. However, the general, nonergodic case is still
   an open
  problem in both the discrete and the continuous settings.''

However, now, the situation has changed. Firstly,  a Tauberian theorem was proposed for differential games in nonergodic case in \cite{Khlopin3}.
  In \cite{Ziliotto}, a very general approach  to proof of Tauberian theorem was proposed for games  with two players with opposite goals
  in  discrete setting.

  In this paper, we show a number of uniform Tauberian theorems for dynamic two-person games with zero sum, in particular, for differential games, for normal form games, and for  games in  discrete setting. Our approach to proofs continues the ideas of  \cite{Lehrer,MN1981,barton}:
   Bellman's optimality principle and the closedness of strategies under concatenation.
 Under these circumstances, the specifics of a game with a saddle point
 (in particular, of a differential game under Isaacs condition) allows to streamline the proof of the
 uniform Tauberian theorem as related to them. Such proof reduces to the two inequalities for the lower and upper games, respectively; following \cite{Bingham} and \cite{Hardy1949}, we call those inequalities One-sided Tauberian theorems.
  After that, the proof of One-sided Tauberian theorem boils down to the application of the suboptimality principle to the strategy that is constructed in a special way through concatenation.

Recall that \cite{barton} provides an example showing that the Tauberian theorem for control problem may not hold if the limits are not uniform over the strongly invariant set of positions. The condition of existence of a uniform limit in the case of Tauberian theorem for dynamic games is certainly as significant. However, this condition may be relaxed for construction of one-sided estimates. Firstly, in construction of such  bounds, one could use any function satisfying the corresponding Dynamic Programming Principle in place of the value function of the original problem. Secondly, for this function, we can replace the condition of existence of the uniform limit with the condition of slow variation (or even monotonicity).


The structure of the paper is as follows.
We start by formulating the Tauberian theorems for normal form games (Theorems~\ref{normal1},\ref{normal2}) in Sect.~\ref{normnormnormalform}.
Then, we consider an axiomatic definition of a rather general game problem statement in the framework of the dynamic model proposed in \cite{barton} (see Sect.~\ref{abstractgame}). At the same place, we formulate Theorems~\ref{maintheoremgeneral},\ref{maintheoremgeneralnew}, the main results of this paper.
 The first of these is based on the existence of a saddle point value in all the considered games, the other one essentially requires the upper value to be greater or equal to the lower value for all of the games.
    In the next Section (Sect.~\ref{side}), we show the connection between the One-sided Tauberian theorem,
    slowly varying functions, and suboptimality principle.
    Sect.~\ref{proof} contains the proofs of   Theorems~\ref{maintheoremgeneral},\ref{maintheoremgeneralnew} and their corollaries, including (the proof of) Theorems~\ref{normal1},\ref{normal2}.
     Sect.~\ref{normalform} is devoted to Tauberian theorem (Theorem~\ref{normtauberian}) for games in discrete setting.
     Most  of the proofs  are located in Appendix. In addition, Appendix~\ref{D} serves to transfer the results obtained for the abstract statement onto the case of differential games (Theorem~\ref{maintheoremdiff}).


  We would like to note that Theorem \ref{maintheoremdiff} and  Theorem~\ref{maintheoremgeneralnew} for more strong assumptions were also proved in \cite{Khlopin3} and \cite{CT15}, respectively.

\section{Normal form game with zero sum}
\label{normnormnormalform}

 Define $\mm{T}\rav\mm{R}_{\geq 0}.$  Assume the following items are given:
 \begin{itemize}
   \item   a nonempty set  $\Omega$;
   \item   a nonempty subset $\mm{K}$ of mappings from $\mm{T}$ to $\Omega$;
   \item   a running cost $g:\Omega\mapsto [0,1];$ for each process $z\in \mm{K},$ assume the  map $t\mapsto g(z(t))$ is  Borel-measurable.
 \end{itemize}

  Let us now  define concatenation on processes.
  Let $\tau\in\mm{T},z',z''\in\mm{K}$ be such that
   $z'(\tau)=z''(0).$ Then, their concatenation $z'\diamond_\tau z''$, a mapping from $\mm{T}$ to $\Omega,$ is defined by the following rule:
\begin{equation}
(z'\diamond_\tau z'') (t)\rav \left\{
 \begin{array} {rcl}        z'(t),       &\mathstrut&     t\leq\tau;\\
                            z''(t-\tau),
                                         &\mathstrut&     t>\tau.
 \end{array}            \right.
\label{3100}
\end{equation}

For all $\omega\in\Omega$, let there be given the non-empty sets $\ct{L}(\omega),\ct{M}(\omega).$
Define sets $\mathfrak{L},\mathfrak{M}$ of  all selectors $\Omega\ni\omega\to l(\omega)\in\ct{L}(\omega),$
$\Omega\ni\omega\to m(\omega)\in\ct{M}(\omega),$
 respectively.
Let, for all   $\omega\in\Omega$, each pair $(l,m)\in\mathfrak{L}\times\mathfrak{M}$ of players' rules
generate a unique process $z[\omega,l,m]\in \mm{K}$  such that $z[\omega,l,m](0)=\omega$.

\begin{theorem}\label{normal1}
For all  $\tau>0$, assume that,  for all $l',l''\in\mathfrak{L}$,
 there exists $l^*\in\mathfrak{L}$
 such that, for all  $\omega\in\Omega$, one has
   \begin{subequations}
\begin{eqnarray}
\label{newa}
\big\{z[\omega,l^*,m']\big|\, m'\in\!\mathfrak{M}\big\}\!=\!
\Big\{
z[\omega,l',m']\diamond_\tau \!z\big[
\omega',l'',m''\big]
\Big|\, m',m''\in\!\mathfrak{M},\omega'\rav z\big[\omega,l',m'\big](\tau)\in\!\Omega\Big\}\subset\mm{K};
\end{eqnarray}
 moreover,  for all $m',m''\in\mathfrak{M}$,
 there exists $m^*\in\mathfrak{M}$
 such that, for all  $\omega\in\Omega$, one has
\begin{eqnarray}
\label{newb}
\big\{z[\omega,l',m^*]\big|\, l'\in\!\mathfrak{L}\big\}\!=\!\Big\{
z[\omega,l',m']\diamond_\tau \!z[\omega',l'',m'']
\Big|\, l',l''\in\!\mathfrak{L},\omega'\rav z\big[\omega,l',m'\big](\tau)\in\!\Omega\Big\}\subset\mm{K}.
\end{eqnarray}
  \end{subequations}
            Assume also that, for each $\lambda,T,h>0$, $\omega\in\Omega$,
    \begin{subequations}
    \begin{eqnarray}
\avg_T(\omega)&\rav&
\sup_{l\in\mathfrak{L}}\inf_{m\in\mathfrak{M}}\frac{1}{T}\int_{0}^T g(z[\omega,l,m](t))\,dt=
\inf_{m\in\mathfrak{M}}\sup_{l\in\mathfrak{L}}\frac{1}{T}\int_{0}^T g(z[\omega,l,m](t))\,dt,\label{2900a}
\\
\bwg_\lambda(\omega)&\rav&
\sup_{l\in\mathfrak{L}}\inf_{m\in\mathfrak{M}}
\lambda\int_{0}^\infty e^{-\lambda t}g\big(z[\omega,l,m](t)\big)\,dt=
\inf_{m\in\mathfrak{M}}\sup_{l\in\mathfrak{L}}
 \lambda\int_{0}^\infty e^{-\lambda t}g\big(z[\omega,l,m](t)\big)\,dt,\label{2900b}\\
\avg_T(\omega)&=&
\sup_{l\in\mathfrak{L}}\inf_{m\in\mathfrak{M}}
 \bigg[\frac{1}{T+h}\int_{0}^h g(z[\omega,l,m](t))\,dt+\frac{T}{T+h}\avg_T\big(z[\omega,l,m](h)\big)\bigg]\label{2950a}\\
&=&
\inf_{m\in\mathfrak{M}}\sup_{l\in\mathfrak{L}}
 \bigg[\frac{1}{T+h}\int_{0}^h g(z[\omega,l,m](t))\,dt+\frac{T}{T+h}\avg_T\big(z[\omega,l,m](h)\big)\bigg],\nonumber\\
\bwg_\lambda(\omega)&=&
\sup_{l\in\mathfrak{L}}\inf_{m\in\mathfrak{M}}
 \bigg[\lambda\int_{0}^h e^{-\lambda t}g\big(z[\omega,l,m](t)\big)\,dt+e^{-\lambda h}\bwg_\lambda\big(z[\omega,l,m](h)\big)\bigg]\label{2950b}\\
&=&
\inf_{m\in\mathfrak{M}}\sup_{l\in\mathfrak{L}}
 \bigg[\lambda\int_{0}^h e^{-\lambda t}g\big(z[\omega,l,m](t)\big)\,dt+e^{-\lambda h}\bwg_\lambda\big(z[\omega,l,m\big](h)\big)\bigg].\nonumber
\end{eqnarray}
\end{subequations}

  Then,
    the following limits exist, are uniform in $\omega\in\Omega$, and coincide
    \begin{equation}\label{6211}
    \lim_{T\uparrow\infty}\avg_T(\omega)
       =\lim_{\lambda\downarrow 0}\bwg_\lambda(\omega)\quad \forall\omega\in\Omega
       \end{equation}
if at least one of these limits exists, and is uniform   in
 $\omega\in\Omega.$
    \end{theorem}
This theorem itself will be proved in Sect.~\ref{proof} as an immediate consequence of Theorem~\ref{maintheoremgeneral}.
In view of that,
conditions  \rref{newa},\rref{newb} express nothing more than the need to test the closedness with respect to concatenation.

 In Theorem~\ref{normal1}, we also need a
 saddle point for all the games considered. This condition can be relaxed, which will also simplify conditions  \rref{newa},\rref{newb}.
\begin{theorem}\label{normal2}
 Assume that, for all $\tau>0$, each of the  sets $\mathfrak{L}$ and $\mathfrak{M}$ is equipped with  a binary operation~$\diamond_\tau$ such that,
 for all $l',l''\in\mathfrak{L},m',m''\in\mathfrak{M}$,
    \begin{subequations}
 \begin{eqnarray}\label{buka}
 z\big[\omega,l'\diamond_\tau l'',m'\diamond_\tau m''\big]=z[\omega,l',m']\diamond_\tau z\big[z[\omega,l',m'](\tau),l'',m''\big]&\quad& \forall \omega\in\Omega,\\
\exists\,l\in\mathfrak{L},m\in\mathfrak{M}\qquad l'\diamond_\tau l=l',m'\diamond_\tau m=m'.&\quad& \label{byka}
 \end{eqnarray}
    \end{subequations}

If there exist the limits in
    \begin{subequations}
\begin{eqnarray}
\lim_{T\uparrow\infty}\sup_{l\in\mathfrak{L}}\inf_{m\in\mathfrak{M}}\frac{1}{T}\int_{0}^T g(z[\omega,l,m](t))\,dt=\lim_{T\uparrow\infty}
\inf_{m\in\mathfrak{M}}\sup_{l\in\mathfrak{L}}\frac{1}{T}\int_{0}^T g(z[\omega,l,m](t))\,dt\quad\forall\omega\in\Omega,\label{naa}
\end{eqnarray}
 in addition, these limits are uniform on $\Omega,$ and coincide, then
all limits in
\begin{eqnarray}
\lim_{\lambda\downarrow 0}\sup_{l\in\mathfrak{L}}\inf_{m\in\mathfrak{M}}\lambda\int_{0}^\infty e^{-\lambda t}g(z[\omega,l,m](t))\,dt=\lim_{\lambda\downarrow 0}
\inf_{m\in\mathfrak{M}}\sup_{l\in\mathfrak{L}}\lambda\int_{0}^\infty
e^{-\lambda t}g(z[\omega,l,m](t))\,dt\ \forall\omega\in\Omega\label{nab}
\end{eqnarray}
    \end{subequations}
exist, are uniform on $\Omega,$ and coincide with the limits in \rref{naa}.

On the other hand, if limits in \rref{nab} exist, are uniform on $\Omega,$ and coincide, then
 the limits in \rref{naa} exist, are uniform on $\Omega,$ and coincide with limits in \rref{nab}.
    \end{theorem}
The theorem will be proved in Sect.~\ref{proof}. For the use of condition \rref{buka} for differential games, refer to \cite[Remark~3.2]{Khlopin3}.

 \section{Abstract dynamic game with zero sum}
\label{abstractgame}

Before exploring the formal definitions, let us sketch a possible interpretation of the necessary formalizations. For stochastic games, a similar statement may be found in \cite{Renault2014}.

 Assume players get some information on state at the current time, but all information on this that is available to the players is contained in a certain signal ~$\omega$; denote the set of all possible signals by~$\Omega$. Since the game develops with time, we can consider the set  $\mm{K},$ which would contain all processes $t\mapsto \omega(t)$ that are possible for the given game. Assume the current value  of  running cost is known at every point of time and, therefore, contained in the signal, i.e., can be described by a function that depends only on~$\omega$.
 Then, by virtue of the known dependence  $t\mapsto \omega(t)$, i.e., in view of the element $z\in\mm{K}$, we can reconstruct the value of the payoff function that realizes.

  Each player also may make private actions according to some rule, as a function of the signal; in this case, there may be some restrictions on the rule's feasibility and on the use of information. We will only assume that the set of such feasible rules is nonempty
  (playable strategy); later, we will also require the existence of $\epsi$-optimal rule for each player. For a fixed initial signal $\omega$, every pair of rules chosen by the players restores some processes   $z\in\mm{K}$ with the property $z(0)=\omega$ (in the case of the normal form games, such a process is unique).
   Denote all such~$z$  by $\Gamma(\omega).$
  On the other hand, each player has the right to publish his decision rule in advance. This act would map to every such rule some subset of ${\mm{K}}$ of the various~$z$ that agree with this rule. Exhausting all the rules, we obtain a set of such subsets, one set for each player
   (${\mA}$ and ${\mB}$ for the first and the second player, respectively).

   If one of the players, let it be the first one, publishes his rule, with this act essentially defining $A\in{\mA}$, let us assume that now the second player has the right to use any information in the choice of initial $\omega$, or rather a process
   $z\in A\cap \Gamma(\omega)$. We (at least, within Theorem~\ref{maintheoremgeneral}) are going to consider the case when mandating the publication for one of the
   players (discrimination against this player)  does not change the value that is guaranteed to him by the  initial information, i.e., when there is a saddle point in such a game.

   Let us now proceed to formal definitions.

 {\bf Dynamic system.}\
 Assume we are given a nonempty set  $\Omega$,  a nonempty subset $\mm{K}$ of mappings from $\mm{T}$ to $\Omega$, and
      a running cost $g:\Omega\mapsto [0,1]$ such that, for each process $z\in \mm{K},$ the  map $t\mapsto g(z(t))$ is  Borel-measurable. For all  $\tau\in\mm{T},z',z''\in\mm{K}$ with
   $z'(\tau)=z''(0),$ we can define the concatenation $z'\diamond_\tau z'':\mm{T}\to\Omega$  by the rule \rref{3100}.

 We would also need the set
    \begin{eqnarray*}
    \Gamma(\omega)\rav\{z\in\mm{K}\,|\,z(0)=\omega\},
     \end{eqnarray*}
     defined for each $\omega\in\Omega$.
This is the set of all feasible processes $z\in\mm{K}$
that begin at $\omega$.

  Let us now  define concatenation on subsets of $\mm{K}$.
  For each pair of  non-empty subsets of $\mm{K}$  and a time $\tau\in\mm{T}$, define
  their concatenation by
\begin{eqnarray}
  A'\diamond_\tau A''&\rav&\{z'\diamond_\tau z''\,|\,z'\in{A}',z''\in A'',z'(\tau)=z''(0)\}\nonumber\\
  &=&\{z'\diamond_\tau z''\,|\,z'\in{A}',z''\in A''\cap \Gamma(z'(\tau))\}. \label{310}
\end{eqnarray}
  To get rid of excessive parentheses, let us hereinafter assume $A\diamond_{\tau'} A'\diamond_{\tau''} A''\rav
  \big(A\diamond_{\tau'} A'\big)\diamond_{\tau''} A''.$

 {\bf Assumptions on strategies.}\
Assume we are given a non-empty family $\mA$
of subsets of the set $\mm{K}$.
   Call a subset $A$ of the set $\mm{K}$ a
  {\it playable strategy} if  we have $A\cap \Gamma(\omega)\neq\pust$  for every initial  $\omega\in\Omega$.

 We hereinafter impose the following conditions on $\mA$:
 \begin{description}
   \item[$(\ct{P})$] $\mA$ is some non-empty set of playable strategies;
   \item[$(\diamond)$] $\mA$ is closed under concatenation $\diamond$: $\forall \tau>0,A',A''\in\mA$
   $A'\diamond_\tau A''\in\mA.$
 \end{description}
Condition $(\ct{P})$ is necessary for all strategies to be applicable for whichever starting information.
Condition~$(\diamond)$  allows the player to switch strategies at some a priori defined time.

{\bf Formalization of lower game.}\
 Consider a two-player lower game with a payoff function $c:\mm{K}\to \mm{R}$. The first player wishes  to maximize
  $c$; the second player wishes to minimize it.
   The first player  also has a family $\mA$ of playable strategies. 

  The lower game is conducted in the following way:
  for a given $\omega\in\Omega,$
  the first player demonstrates a set  $A\in\mA$, and
  then the second player chooses a process ${z\in A}\cap\Gamma(\omega)$.
  The value function of this game is
  \begin{subequations}
  \begin{eqnarray}
 \cym[c](\omega)\rav\sup_{A\in\mA}\inf_{z\in A\cap \Gamma(\omega)}c(z)\qquad \forall\omega\in\Omega.
 \label{280}
\end{eqnarray}
 Note that  for every playable strategy $A$, $A\cap \Gamma(\omega)\neq\pust;$ therefore, this definition is valid
 if $c$ is bounded.

    \begin{definition}
    For each positive $\varepsilon,$
   let us say that the lower game of \rref{280}  has {\it an $\varepsilon$-optimal strategy} $A\in\mA$ if
     $\inf_{z\in A} c(z)\geq\cym[c](z(0))-\varepsilon.$
 \end{definition}

{\bf Formalization of upper game.}\
  We still have two players, the first player  maximizes
  the payoff~$c$, whereas the second player minimizes it.
  Let the second player
    also have a family $\mB$ of playable strategies.

  The upper game is conducted in the following way.
  Given $\omega\in\Omega,$
  let the second demonstrate some   set $B\in\mB$,
  then let the first player choose a process ${z\in B}\cap\Gamma(\omega)$.
  The value function of the upper game is
  \begin{eqnarray}
 \cyp[c](\omega)\rav\inf_{B\in\mB}\sup_{z\in B\cap \Gamma(\omega)}c(z)\qquad \forall\omega\in\Omega.
 \label{580}
\end{eqnarray}
 \end{subequations}
    \begin{definition}
    For each positive $\varepsilon,$
   let us say that the upper game of \rref{580}  has {\it an $\varepsilon$-optimal strategy} $B\in\mB$ if
     $\sup_{z\in A} c(z)\leq\cyp[c](z(0))+\varepsilon.$
 \end{definition}

Note that
$1-\cyp[c](\omega)=\sup_{B\in\mB}\inf_{z\in B\cap \Gamma(\omega)}\big(1-c(z)\big)$  for all $\omega\in\Omega,$
 i.e., the upper game with payoff $c$ and with set $\mB$ of the second player's strategies
 differs from
 the lower game with payoff $1-c$ and with set $\mB$ of the first player's strategies only in its sign. Consequently, all definitions and statements below will mostly be given for the
  lower game. For the upper game families, they can be obtained by the replacement
  $g^-\rav 1-g,$ $\mm{A}^-\rav \mm{B}.$

In what follows, fix $\mA$ and $\mB$;
  to every
 payoff $c$, it is possible to map the pair of corresponding games by rules  \rref{280} and \rref{580}.

  \begin{definition}
    For 
    a payoff $c:\mm{K}\to\mm{R}$,
     let us say that the corresponding  games   have {\it a  saddle point} if $\cym[c](\omega)=\cyp[c](\omega)$ for all $\omega\in\Omega.$
    \end{definition}
{\bf On various payoffs.}\
 Let us now define time average $\av_T(z)$ and discount average $\bw_\lambda(z)$ for each process $z\in \mm{K}$
 by the rules:
 \begin{eqnarray}
    \label{248}
 \av_T(z)\rav\frac{1}{T}\int_{0}^T g(z(t))\,dt,\quad
 \bw_\lambda(z)\rav\lambda\int_{0}^\infty e^{-\lambda t} g(z(t))\,dt\qquad \forall z\in\mm{K},T>0,\lambda>0.
 \end{eqnarray}
   Note that the definitions are valid, and  the means lie within  $[0,1].$
 The functions $\av_T,\bw_\lambda$ will be treated as the payoff functions in their respective games.
 In particular, for any $T,\lambda>0$, we obtain the values
\begin{subequations}
  \begin{eqnarray}
 \avm_T(\omega)\!\!\!\!&\rav&\!\!\!\!\cym[\av_T](\omega)\!=\!\sup_{A\in\mA}\inf_{z\in A\cap \Gamma(\omega)}\!\av_T(z), \avp_T(\omega)\!\!\rav\!\!\cyp[\av_T](\omega)\!=\!\inf_{B\in\mB}\sup_{z\in B\cap \Gamma(\omega)}\!\av_T(z),\;\forall\omega\in\Omega,
 \label{280a}\\
 \bwm_\lambda(\omega)\!\!\!\!&\rav&\!\!\!\!\cym[\bw_\lambda](\omega)\!=\!\sup_{A\in\mA}\inf_{z\in A\cap \Gamma(\omega)}\!\bw_\lambda(z),\
 \!\!\bwp_\lambda(\omega)\!\rav\!\cyp[\bw_\lambda](\omega)\!=\!\!\inf_{B\in\mB}\!\sup_{z\in B\cap \Gamma(\omega)}\!\bw_\lambda(z),   \forall\omega\in\Omega. \label{280b}
\end{eqnarray}
\end{subequations}
We will also need Bolza-type payoffs.
For each map  $U:\mm{R}_{>0}\times\Omega\to\mm{R}$ and positive $h,T,\lambda$, define the following  payoffs:
 \begin{subequations}
\begin{eqnarray}
\label{533}
 \hat{c}^U_{h,T}(z)&\rav& 
 \frac{1}{T+h}\int_{0}^h g(z(t))\,dt+\frac{T}{T+h}U_T(z(h))\qquad\forall  z\in\mm{K};\\
\label{633}
\check{c}^U_{h,\lambda}(z)&\rav& 
\lambda\int_{0}^h e^{-\lambda t}g(z(t))\,dt+e^{-\lambda h}U_{\lambda}(z(h))\qquad\forall z\in\mm{K}.
\end{eqnarray}
 \end{subequations}
  For the following theorem, in \rref{533} and \rref{633}, we need to use $U=\avm=\avp,U=\bwm=\bwp$,
 respectively.

{\bf Uniform Tauberian theorem for games with a saddle point.}\
    \begin{theorem}
     \label{maintheoremgeneral}
            Let $\mA$ and $\mB$ satisfy conditions $(\ct{P})$,$(\diamond)$.

            Assume that for each $\lambda,T,h>0$, for each of the following payoffs
             $\av_T$, $\hat{c}^{\avm}_{h,T},$
             $\bw_\lambda$, $\check{c}^{\bwm}_{h,\lambda}$
             the corresponding  games have saddle points and
$\epsi$-optimal player's strategies from $\mA,\mB$ respectively for all $\epsi>0$.

  Then,   the following limits exist, are uniform in $\omega\in\Omega$, and coincide
    \begin{equation}\label{621}
    \lim_{T\uparrow\infty}\avp_T(\omega)=\lim_{T\uparrow\infty}\avm_T(\omega)
       =\lim_{\lambda\downarrow 0}\bwp_\lambda(\omega)=\lim_{\lambda\downarrow 0}\bwm_\lambda(\omega)\quad \forall\omega\in\Omega
       \end{equation}
if at least one of these limits exists, and is uniform   in
 $\omega\in\Omega.$
    \end{theorem}
For the proof of this theorem, refer to Sect.~\ref{proof}.
The condition of existence of a
saddle point can be relaxed, see Corollary~\ref{maintheoremgeneral1},
also in Sect.~\ref{proof}.

The conditions of Theorem~\ref{maintheoremgeneral} mostly deal with the value functions, however, a similar theorem can be formulated in terms of  the capabilities of players.
To this end, consider the following:

{\bf  Additional conditions for capabilities' sets.}

  In addition to  concatenation, let us also define time-shift.
  For a time $\tau\in\mm{T}$ and a process $z\in\mm{K}$, define the function
   $z_\tau:\mm{T}\mapsto\Omega$ by the following rule:
\begin{equation}
\label{296}
z_\tau(t)=z(t+\tau) \qquad\forall t\in\mm{T}.
\end{equation}

  We say that a family $\mA$ of playable strategies
 allows the separation of $\omega$ (at the initial time) if, for each mapping
   $\eta:\Omega\to \mA$,
 \begin{equation*}
    \exists A\in\mA\ \forall \omega\in\Omega\qquad \eta(\omega)\cap\Gamma(\omega)=A\cap\Gamma(\omega).
 \end{equation*}

   Call  a family $\mA$ closed under backward shift if,
    for all $A\in\mA$, $\tau>0$, there exists a $A'\in\mA$ such that
  $A=A\diamond_\tau A'.$ It is easy see that $A=A\diamond_\tau A'$ implies $\{z_\tau\,|\,z\in A\}\subset A'.$

  We say that families $\mA$ and $\mB$  are compatible  if,
  for all $\omega\in\Omega$,$A\in\mm{A},B\in\mm{B},$
  $A\cap B\cap \Gamma(\omega)$  is non-empty.

So, let us also  define the following conditions:
 \begin{description}
   \item[$(\ct{C})$] $\mA$ and $\mB$   are compatible;
   \item[$(\omega)$] $\mA$ allows the separation of $\omega$ (at the initial time);
   \item[$(\tau)$] $\mA$   is closed under backward shift.
 \end{description}

$(\ct{C})$ is a stronger variant of condition $(\ct{P})$; the former is relatively often used for the games with a
saddle point, see, for example \cite[Subsect.~VIII.3]{Bardi}. We will use  it in Lemma~\ref{firstlemm2}.
Condition~$(\omega)$ lets the player use the data on the initial $\omega$;
under this condition, we can always provide $\epsi$-optimal players' strategies (see Lemma~\ref{firstlemm1}).
$(\diamond)$\&$(\omega)$ let the player plan the switch of strategies at a given time in advance, at the beginning of the game; the switch is based on the information about the state that would realize at the time of the switch. In particular, this switching can be used to construct, in the framework of the Tauberian theorem, the near-optimal strategies, see Remark~\ref{555} and \rref{UT},\rref{Ul}.
Condition $(\tau)$ means that an action is admissible for this information
 at a positive time if this action is admissible for this information  at zero time; below, we will only use this condition to prove Bellman's optimality  principle
(see  Lemma~\ref{crucial}).

{\bf Uniform Tauberian theorems for games  without a saddle point.
}\
    \begin{theorem}
     \label{maintheoremgeneralnew}
            Let $\mA$ and $\mB$ satisfy conditions $(\ct{C})$,$(\diamond)$,$(\omega)$,$(\tau)$.

  If,  either for lower and upper games with payoffs  $\av_T$ $(T>0)$, or for lower and upper games with payoffs $\bw_\lambda$ $(\lambda>0)$,  limits of their values  (in \rref{621}) exist, are uniform in $\omega\in\Omega$, and coincide,
  then,   all limits in \rref{621} exist, are uniform in $\omega\in\Omega$, and coincide.
    \end{theorem}
Our main objective for two next Sections is to prove these theorems.

 \section{One-sided Tauberian theorems.}
\label{side}

{\bf On game families.}\
 Below, we will consider various directed sets
  of payoffs, which we will index by positive numbers.
  For example,
  $\bw_\lambda$ $(\lambda\downarrow 0),$
 $\av_T$ $(T\uparrow\infty),$
  $\bw_\lambda$ $(\lambda\downarrow 0),$
 $\av_T$ $(T\uparrow\infty),$
 $\hat{c}^U_{h,T}$ $(T\uparrow\infty),$ $\check{c}^U_{h,\lambda} (\lambda\downarrow 0)$
 for some $h>0,$ $U:\mm{R}_{>0}\times\Omega\to \mm{R}.$
 For brevity, for all of these payoff families, we might also use the notation
 $\nu_\gamma$ $(\gamma\to \gamma_*),$ 
 with $\gamma>0,$ $\gamma_*\in\{0,+\infty\},$ and
 payoffs $\nu_\gamma.$
 For all $\gamma>0$, the set $[\gamma,\infty)$ (respectively, $(0,\gamma]$)
 is called a neighborhood of $\gamma_*$ if  $\gamma_*=+\infty$
 (if $\gamma_*=0$).

  {\bf On  guaranties and suboptimality principles.}\
  Consider functions $S:\Omega\to \mm{R}$, $U:\mm{R}_{>0}\times\Omega\to \mm{R}$.
    \begin{definition}
  Let us say that a lower game
  \rref{280} (respectively, or \rref{280a}, or \rref{280b}) {\it has  a guarantee} $S$ if
  $\cym[c](\omega)\geq S(\omega)$
  (respectively,
   or $\avm_T(\omega)\geq S(\omega)$, or $\bwm_\lambda(\omega)\geq S(\omega)$) for all $\omega\in\Omega$.

   Let us say that a  guarantee $S$ of lower game \rref{280} with payoff $c$
   {\it is protected}  if there exists a strategy $A\in\mA$
  satisfying
  $c(z)\geq S(z(0))$ for all $z\in A.$

 Let us say that a lower game family with payoffs $\nu_\gamma(\gamma\to\gamma_*)$
 has {\it an  asymptotic guarantee} $U$ ({\it a protected asymptotic guarantee}) if, for each $\epsi>0$, there exists a
  neighborhood of $\gamma_*$
 such that $U_\gamma-\epsi$, as function from $\Omega$ to $\mm{R}$, is a guarantee (a protected guarantee)
  for lower game  with  payoff $\nu_\gamma$ for all $\gamma$ from this neighborhood.
   \end{definition}
Note that each function $S:\Omega\to \mm{R}$ may be regarded as a function
from $\mm{R}_{>0}\times\Omega$ to $\mm{R}$, which does not depend on the first argument; therefore, $S$
may be
an  asymptotic guarantee.

\begin{subequations}
    \begin{definition}
 Let us say that a function $U:\mm{R}_{>0}\times\Omega\to\mm{R}$
  is
  {\it a  subsolution} of the lower game family
  with payoffs $\av_T(T\uparrow\infty)$
 if, for every $\varepsilon>0$,
  there exists a positive  $\bar{T}$ such that, for all $h>0,T>\bar{T}$,
 there exists a strategy $A\in\mA$ such that
\begin{eqnarray}
 U_{T+h}(z(0))-\varepsilon\leq
  \hat{c}^U_{h,T}(z)=\frac{1}{T+h}\int_{0}^{h}g(z(t))\,dt+\frac{T}{T+h}U_{T}(z(h))
   \qquad \forall z\in A,
 \label{sol207}
\end{eqnarray}
   i.e.,\,\,for all $h\!>\!0,T\!>\!\bar{T}$,
  the function $U_{T\!+\!h}\!-\!\varepsilon$ is a protected guarantee of the lower game with payoff\,\,$\hat{c}^U_{h,T}$.

 Let us say that a function $U:\mm{R}_{>0}\times\Omega\to\mm{R}$
  is
  {\it a subsolution} of the lower game family
  with payoffs $\bw_\lambda(\lambda\downarrow 0)$
  if, for every $\varepsilon>0$,
  there exists a positive  $\bar\lambda$ such that, for all positive $\lambda<\bar\lambda,$ for all $h>0$,
 there exists a strategy  $A\in\mA$ such that
\begin{eqnarray}
 U_{\lambda}(z(0))-\varepsilon\leq\check{c}^U_{h,\lambda}(z)=
  \lambda\int_{0}^{h}e^{-\lambda t}g(z(t))\,dt+e^{-\lambda h} U_{\lambda}(z(h))
 \qquad \forall z\in A,
 \label{sol307}
\end{eqnarray}
   i.e., for all positive $\lambda<\bar\lambda,$ for all $h>0$,
 $U_{\lambda}-\varepsilon$ is a protected guarantee of lower game with payoff $\check{c}^U_{h,\lambda}.$
    \end{definition}
For a similar definition, refer to \cite[Sect.~VI.4]{Bardi},
 suboptimality principle \cite[Definition III.2.31]{Bardi} (also referred to as `stability with respect to second player'~\cite{ks}, and \cite[Sect.~VI.4]{Bardi} for discrete problems).
\end{subequations}
\begin{remark}\label{888}
 Let a payoff family $\nu_\gamma$ $(\gamma\to \gamma_*)$ be either
 $\av_T$ $(T\uparrow \infty)$, or   $\bw_\lambda$ $(\lambda\downarrow 0)$,
 Assume $U,U'$ from $\mm{R}\times\Omega$ to $\mm{R}$ has a common limit as $\gamma\to\gamma_*,$
 and this limit is uniform on $\Omega$.
 Then,
\begin{enumerate}
  \item $U$ is  an  asymptotic guarantee of the lower game family with payoffs $\nu_\gamma(\gamma\to\gamma_*)$ iff
 $U'$ is the same;
  \item $U$ is  a protected asymptotic guarantee of this  lower game family
   iff
 $U'$ is the same;
  \item $U$ is  a subsolution of this lower game family
 iff
 $U'$ is the same.
\end{enumerate}
\end{remark}

{\bf One-sided Tauberian theorems for lower game families.}\

    In Appendix \ref{A}, we prove  the following proposition:
 \begin{subequations}
     \begin{proposition}
     \label{av_bw}
       Let $\mA$  satisfy conditions $(\ct{P})$,$(\diamond)$.

       Let a bounded from above function $U:\mm{R}_{>0}\times\Omega\to\mm{R}$
        be a subsolution for the lower game family with
        payoffs $\av_T (T\uparrow\infty)$, in particular, let
                 \rref{sol207} hold.

      Let  also $U$ 
      satisfy
 \begin{equation}
 \label{slowlyT}
 \liminf_{T\uparrow\infty} \inf_{\omega\in\Omega} \big(U_{pT}(\omega)-U_{T}(\omega)\big)\geq 0\qquad\forall p>1.
 \end{equation}

      Then,
      for every $\varepsilon>0$, there exists a natural $N$ such that, for all positive $\lambda<1/N,$ a function  $U_{1/\lambda}-\varepsilon$ is a protected guarantee of the lower game with payoff $\bw_\lambda;$
      in particular,
       \begin{equation*}
       \bwm_{\lambda}(\omega)\geq U_{1/\lambda}(\omega)-\varepsilon\qquad\forall \omega\in\Omega,\lambda\in(0,1/N).
       \end{equation*}
     \end{proposition}

       In Appendix \ref{B}, we prove  the following proposition:
 \begin{proposition}
 \label{bw_av}
       Let $\mA$  satisfy conditions $(\ct{P})$,$(\diamond)$.

       Let a bounded from above function $U:\mm{R}_{>0}\times\Omega\to\mm{R}$
        be
        a subsolution for the lower game family
         with payoffs $\bw_\lambda (\lambda\downarrow 0),$ in particular,
let \rref{sol307} hold.

      Let  also $U$  satisfy
  \begin{equation}
 \label{slowlyl}
 \liminf_{\lambda\downarrow 0} \inf_{\omega\in\Omega}
 \big(U_{\lambda}(\omega)-U_{p\lambda}(\omega)\big)\geq 0\qquad\forall p>1.
 \end{equation}

       Then, for every
       $\varepsilon>0$,  there exists a natural $N$ such that, for all positive $T>N,$
       a function $U_{1/T}-\varepsilon$ is a protected guarantee of the lower game with payoff $\av_T;$
      in particular,
       $$\avm_T(\omega)\geq U_{1/T}(\omega)-\varepsilon\qquad  \forall \omega\in\Omega,T>N.$$
     \end{proposition}
Note that \rref{slowlyl} holds for $U$ if
 either $U$ is decreasing in $\lambda$ for all $\omega\in\Omega$,
 or $U$ has  limit as $\lambda\downarrow 0$, and this limit is uniform in  $\omega\in\Omega,$
 or $U$ is a sum of such functions.
Analogously, \rref{slowlyT} holds for $U$ if
 either $U$ is increasing in $T$ for all $\omega\in\Omega$
 or $U$ has  limit as $T\uparrow \infty$, and this limit is uniform in  $\omega\in\Omega,$
 or $U$ is a sum of such functions.

The inequalities similar to \rref{slowlyT} will be found in  Uniform Tauberian Theorem,
(see definitions of  slowly decreasing, slowly varying, for example, in \cite[Definition 4.1.4]{Bingham},\cite[Sect.6.2]{Hardy1949}).
For stochastic  games the similar inequalities can see in
 \cite{MN1981} (Condition~$(3^*)$ from Theorem 4.1).

Note that, in addition to uniform and exponential payoff families, the Tauberian theorems can be formulated for   arbitrary
probability densities. The corresponding results  are seen  for discrete time systems \cite{1993},\cite{Renault2013}, for optimal control \cite{LQR},
for games \cite{CT15}. In \cite{CT15}, it is made based on the corresponding   One-sided Tauberian theorems similar to Propositions~\ref{av_bw},\ref{bw_av}.
 \end{subequations}

    \begin{corollary}
     \label{maincorollary}
            Let $\mA$  satisfy conditions $(\ct{P})$,$(\diamond)$.
 Let the functions $\avm$, $\bwm$
       be subsolutions for the lower game families with payoffs $\av_T(T\uparrow\infty)$
       and with payoffs $\bw_\lambda(\lambda\downarrow 0)$, respectively.

       If there exists
    a limit of $\avm_T$ that is uniform  on $\Omega$:
\begin{equation*}
S_*(\omega)=\lim_{T\uparrow\infty}\avm_T(\omega)\qquad \forall \omega\in\Omega,
\end{equation*}
       then, this limit $S_*$ is a common protected asymptotic guarantee  for
       lower game families  with  payoffs $\av_T(T\uparrow\infty)$ and  $\bw_\lambda(\lambda\downarrow 0)$.

       If there exists
 a limit of $\bwm_\lambda$ that is uniform  on $\Omega$:
\begin{equation*}
S_*(\omega)=\lim_{\lambda\downarrow 0}\bwm_\lambda(\omega)\qquad \forall \omega\in\Omega,
\end{equation*}
        then,  this limit $S_*$ is
      a  common protected asymptotic guarantee  for
       lower game families with  payoffs $\av_T(T\uparrow\infty)$ and  $\bw_\lambda(\lambda\downarrow 0)$.
    \end{corollary}
\proof of Corollary~\ref{maincorollary}.
       Assume that there exists
    a limit of $\avm_T$  that is uniform  on $\Omega$. The uniformity of this limit implies \rref{slowlyT} for $U\equiv\avm,$ also, we see that $\avm, S_*$ are asymptotic guaranties  for the
       lower game family with payoffs $\av_T(T\uparrow\infty).$

       We claim that $S_*$ is a protected  asymptotic guarantee  for the
       lower game family with $\bw_\lambda.$
        By condition, $\avm$ is a subsolution
       for this
       lower game family. Therefore, $S_*$ is the same by Remark~\ref{888}.
     Now, by Proposition~\ref{av_bw},
 $S_*$ is a protected asymptotic guarantee  for
       lower game family with payoffs $\bw_\lambda.$

      The second part is  proved analogously by Proposition~\ref{bw_av},
      i.e.,  since $S_*$ is an asymptotic guarantee  for the
       lower game family with payoffs  $\bw_\lambda,$
        $S_*$ is  also  a protected asymptotic guarantee  for the
       lower game family with $\av_T.$ By the first part of the proof, now
         $S_*$ is a protected asymptotic guarantee  for
       the lower game family with payoffs  $\bw_\lambda.$
       The proof of the same fact for the payoff family $\av_T$ $(T\uparrow\infty)$ is similar.
      \bo

\section{Proofs of main results. Their corollaries.}
\label{proof}

{\bf Auxiliary lemmas. The formulations.}

Similarly  \cite[Remark III.2.7]{Bardi}, for  continuous dynamic systems with payoff $\bw_\lambda$,
we obtain:
\begin{lemma}
     \label{lastlemma}
            Let the set $\mA$ satisfy conditions $(\ct{P})$,$(\diamond)$.
 Let $\epsi,T,\lambda$ be positive.
 Then,

   1)  $\cym\big[\hat{c}^{U}_{h,T}\big]-\epsi$ is a guarantee of lower game with payoff $\av_{T+h}$ for all $h>0$ if
     $U_T-\epsi$ is a protected guarantee of lower game with payoff $\av_T;$

    2)
    $\cym\big[\check{c}^{U}_{h,\lambda}\big]-\epsi$ is a guarantee of lower game with payoff $\bw_\lambda$ for all $h>0$
    if
         $U_\lambda-\epsi$ is  a protected guarantee of lower game with payoff $\bw_\lambda.$
\end{lemma}
 \begin{lemma}
 \label{crucial}
       Let the set $\mA$
       satisfy conditions  $(\ct{P})$,$(\tau)$. Then,

       1)  $\avm$
       is a  subsolution for the lower game family
       with payoffs $\av_T(T\uparrow\infty)$
       if $\avm$ is a protected asymptotic guarantee for
       this lower game family;

       2)     $\bwm$
       is a   subsolution for the lower game family with payoffs $\bw_\lambda(\lambda\downarrow 0)$
       if $\bwm$ is a protected asymptotic guarantee for
       this lower game family.
\end{lemma}
We will also need the two lemmas that are essentially contained in  \cite[Lemma~2.1]{CT15}.
\begin{lemma}
 \label{firstlemm1}
  Let the set  $\mA$ satisfy conditions $(\ct{P}),(\omega)$. Let a mapping $c:\mm{K}\to \mm{R}$ be bounded.
  Then, for all $\epsi>0$, $\cym[c]-\epsi$ is a protected guarantee for lower game with payoff $c$.
\end{lemma}
\begin{lemma}
 \label{firstlemm2}
  Let the sets $\mA,\mB$ satisfy condition $(\ct{C})$.
  Then, $\cym[c]\leq\cyp[c]$ for a bounded payoff $c:\mm{K}\to \mm{R}$.
\end{lemma}
The  proofs of all lemmas are located in Appendix~\ref{AA}.

  Note that, for a payoff family $\nu_\gamma$ $(\gamma\to \gamma_*)$, if lower and upper game families with $\nu_\gamma(\gamma\to\gamma_*)$ have a~common asymptotic guarantee $S_*:\Omega\to \mm{R},$ then,
  for all $\epsi>0$,  for all $\gamma$ from some neighborhood of~$\gamma_*$,
  we obtain
 \begin{eqnarray}\label{889}
 \cyp[\nu_\gamma](\omega)-\epsi\leq S_*(\omega)\leq \cym[\nu_\gamma](\omega)+\epsi\qquad\forall\omega\in\Omega. \end{eqnarray}
 Conceptually, the proofs of Theorems~\ref{maintheoremgeneral},\ref{maintheoremgeneralnew} and their corollaries consist of the two parts. By  checking the conditions of one of  Propositions \ref{av_bw} and \ref{bw_av}, we prove
 \rref{889}. Its converse inequality, in the case of Theorem~\ref{maintheoremgeneralnew}, is based on  Lemma~\ref{firstlemm2} and, in the case of Theorem~\ref{maintheoremgeneral}, on the existence of a
 saddle point for games with payoffs $\av_T,\bw_\lambda.$

{\bf Proof of Theorem~\ref{maintheoremgeneral}}

By condition,  $\avm_T\equiv\avp_T$, $\bwm_\lambda\equiv\bwm_\lambda$ for all $\lambda,T>0$.
 Therefore, for each $h,\lambda,T>0$, the games with the payoff  $\hat{c}^{\avm}_{h,T}$ and the games with the payoffs $\hat{c}^{\avp}_{h,T}$
 coincide completely. The same is true for the pair of payoffs
 $\check{c}^{\avm}_{h,T}$, $\check{c}^{\avp}_{h,T}.$
  Thus,  we can consider
  the conditions of the theorem with respect to lower and upper games to be totally symmetric.

If at least one of limits from \rref{621} exists and is uniform   on $\Omega,$
then,
    either the limit of $\avm_T=\avp_T$ as $T\uparrow\infty,$
or the limit of $\bwm_\lambda=\bwp_\lambda$ as $\lambda\downarrow 0$
 exists, and
        is uniform  in $\omega\in\Omega.$
 Assume that it is the limit $S_*$ of $\avm_T$, $\avp_T$ as $T\uparrow\infty.$


 We claim that $S_*$ is a protected asymptotic guarantee
for the lower game family with $\bw_\lambda(\lambda\downarrow 0)$.

Fix a $\epsi, T>0.$
  By condition,  we can  find $\epsi$-optimal strategies for payoff $\av_T.$
  So, $\avp_{T}+\epsi$ is a protected guarantee of upper game with payoff $\av_T.$
Using Lemma~\ref{lastlemma} for upper games, we obtain
 $$
 \avp_{T+h}\leq \cyp\big[\hat{c}^{\avp}_{h,T}\big]+\epsi$$
 for all $h>0$.
 By hypothesis, for each $h>0$, there exists $\epsi$-optimal strategy $A^h\in\mA$ for payoff
$\hat{c}^{\avp}_{h,T};$ now,
$$\avp_{T+h}(z(0))-2\epsi\leq\cym\big[\hat{c}^{\avp}_{h,T}\big](z(0))-\epsi\leq
\hat{c}^{\avp}_{h,T}(z)
   \qquad \forall z\in A^h.
  $$
So, \rref{sol207}   with $U=\avp$ holds for all positive $h.$
Thus, $\avp$ is  a subsolution of lower  game family with payoffs $\av_T$.

Then, by Proposition \ref{av_bw}, $\avp$
 is a protected asymptotic guarantee
for lower game family with $\bw_\lambda(\lambda\downarrow 0)$.
Now, by Remark~\ref{888}, $S_*$ is the same.

Analogously, by  symmetry with respect to lower and upper games,
one can prove that $S_*$, as the uniform limit of $\avp$,  is a protected asymptotic guarantee
for upper game family with $\bw_\lambda(\lambda\downarrow 0)$.

 So,   $S_*$ is an  asymptotic guarantee
for lower and upper game families with $\bw_\lambda(\lambda\downarrow 0)$, i.e.,
\rref{889} holds for this family. Now, by $\bwp\equiv\bwm$,  $S_*$ is the limit of  $\bwm_\lambda$ and $\bwp_\lambda$, and this limit is uniform.

  The second part of the proof is very similar, it is only necessary to swap  $\lambda$ and $T$, $\av_T$ and $\bw_\lambda$ everywhere, and use
  $\check{c}_{h,\lambda},$   Proposition~\ref{bw_av} instead of $\hat{c}_{h,T}$,
  Proposition~\ref{av_bw}.
  \bo


    Let us make Theorem~\ref{maintheoremgeneral} more precise. Recall that $\gamma_*\in\{0+,+\infty\}.$
        \begin{definition}
    For   payoff family  $\nu_\gamma$ $(\gamma\to \gamma_*)$,
     let us say that this family
       has {\it an asymptotic saddle point} if  the limit of $|\cym[\nu_\gamma](\omega)-\cyp[\nu_\gamma](\omega)|$ as $\gamma\to \gamma_*$
       exists, is equal to $0$, and
        is uniform  for $\omega\in\Omega$.%
    \end{definition}
    \begin{definition}
    For a monotonic
     function $\varkappa:\mm{R}_{>0}\to\mm{R}_{>0}$ with   $\varkappa(\gamma)\to 0$ as $\gamma\to\gamma_*,$
    let us say that $\varkappa$ is {\it a precision}
    of family of payoffs $\nu_\gamma$ $(\gamma\to \gamma_*)$ if
     the lower and upper games with the payoff $\nu_\gamma$ have
     $\varkappa({\gamma})$-optimal strategies for each player.
 Let us say that $\varkappa$ is {\it a common precision}
    for a set of payoff families if
    $\varkappa$ is a  precision for each  payoff family from this set.
    \end{definition}
    \begin{remark}\label{777}
    Consider a family of payoffs $\nu_\gamma$ $(\gamma\to \gamma_*)$.
     Any monotonic function $\varkappa:\mm{R}_{>0}\to\mm{R}_{>0}$ with   $\varkappa(\gamma)\to 0$ as $\gamma\to\gamma_*$
    is actually  a precision
    of this family
     if all corresponding games have $\epsi$-optimal players' strategies for a positive $\epsi.$
   \end{remark}
    \begin{corollary}
     \label{maintheoremgeneral1}
            Let $\mA$ and $\mB$ satisfy conditions $(\ct{P})$,$(\diamond)$.
            Assume that both payoff families  $\av_T(T\uparrow \infty)$,
             $\bw_\lambda (\lambda\downarrow 0)$
              have  a precision and an asymptotic saddle point.
          Also, for each of set of payoff families
         $$\{\hat{c}^{\avm}_{h,T}(T\uparrow \infty)\,|\,h>0\},\quad
         \{\hat{c}^{\avp}_{h,T}(T\uparrow \infty)\,|\,h>0\},\quad
         \{\check{c}^{\bwp}_{h,\lambda}(\lambda\downarrow 0)\,|\,h>0\},\quad
         \{\check{c}^{\bwp}_{h,\lambda}(\lambda\downarrow 0)\,|\,h>0\},$$
         let there exist a common precision for this set.

   Then,   all limits in \rref{621} exist, are uniform in $\omega\in\Omega$, and coincide
   if at least one of these limits exists and is uniform   in
 $\omega\in\Omega.$
    \end{corollary}
    \proof of Corollary~\ref{maintheoremgeneral1}.
By condition, the limit of $|\avm_T(\omega)-\avp_T(\omega)|$ as $T\uparrow\infty$
and the limit of $|\bwm_\lambda(\omega)-\bwm_\lambda(\omega)|$ as $\lambda\downarrow 0$
equal $0$ and are uniform in $\omega\in\Omega.$

Now,  if at least one of limits in \rref{621} exists, and is uniform   on $\Omega,$
then,
    either the limit of $\avm_T$ and $\avp_T$ as $T\uparrow\infty,$
or the limit of $\bwm_\lambda$ and $\bwp_\lambda$ as $\lambda\downarrow 0$,
 exists and
        is uniform  in $\omega\in\Omega.$
 Assume that it is the limit $S_*$ of $\avm_T$ and $\avp_T$ as $T\uparrow\infty.$

  By condition,  there exists a precision function  $\hat{\varkappa}:\mm{R}_{>0}\to\mm{R}_{>0}$ with $\hat{\varkappa}(T)\to 0$ as $T\uparrow\infty$ such that,
  $|\avm_T(\omega)-\avp_T(\omega)|<\hat{\varkappa}(T)$ for all $\omega\in\Omega.$
  As corollary,  $|\hat{c}^{\avm}_{h,T}(z)-\hat{c}^{\avp}_{h,T}(z)|<\hat{\varkappa}(T)$ for all $z\in\mm{K},h>0.$
  In addition, by condition, increasing the function $\hat{\varkappa}$ if necessary (provided that $\hat{\varkappa}(+\infty)=0)$,  we can propose that, for all $T,h>0$,
  games with payoffs $\hat{c}^{\avp}_{h,T},\hat{c}^{\avp}_{h,T}$ also have $\hat{\varkappa}(T)$-optimal  players' strategies.
  Now, we can consider the conditions of the theorem  to be totally symmetric with respect to lower and upper games.

 We claim that $S_*$ is a protected asymptotic guarantee
for  lower game family with $\bw_\lambda(\lambda\downarrow 0)$.
 Fix a positive $T.$
  By condition,  we can  find $\hat{\varkappa}(T)$-optimal strategies for the payoff $\av_T.$
Using Lemma~\ref{lastlemma} for upper games, we obtain
 $$
 \avp_{T+h}-\hat\varkappa(T)\leq \cyp\big[\hat{c}^{\avp}_{h,T}\big]\qquad  \forall h>0.
$$
 By the choice of $\hat\varkappa$, for all $h>0$, there exists $\hat\varkappa(T)$-optimal strategy $A^h\in\mA$ for payoff
$\hat{c}^{\avp}_{h,T};$ now
$$\avp_{T+h}(z(0))-2\hat{\varkappa}(T)\leq\cym\big[\hat{c}^{\avp}_{h,T}\big](z(0))-\hat\varkappa(T)\leq
\hat{c}^{\avp}_{h,T}(z)
   \qquad \forall z\in A^h.
  $$
So, \rref{sol207}   with $U=\avp$ hold for all positive $h.$
Thus, $\avp$ is  a subsolution of lower game  family with payoffs $\av_T(T\uparrow \infty)$.
Now, by Remark~\ref{888}, $S_*$ is the same.
By Proposition \ref{av_bw}, $S_*$
 is a protected asymptotic guarantee
for lower game family with payoffs $\bw_\lambda(\lambda\downarrow 0)$.

Analogously, by symmetry with respect to lower and upper games,
one can prove that $S_*$, as the uniform limit of $\avp$,  is a protected asymptotic guarantee
for upper game family with $\bw_\lambda(\lambda\downarrow 0)$.

 Now,   $S_*$ is an  asymptotic guarantee
for lower and upper game families with $\bw_\lambda(\lambda\downarrow 0)$, i.e.,
\rref{889} holds for this family. Since $|\bwp_\lambda(\omega)-\bwm_\lambda(\omega)|$ as $\lambda\downarrow 0$
uniformly (on $\Omega$) tends to $0$,    $S_*$ is the limit of  $\bwm_\lambda$ and $\bwp_\lambda$, and this limit is uniform.

  The second part of the proof is analogous, it suffices to swap $\lambda$ and $T$, $\av_T$ and $\bw_\lambda$,  $\lambda\downarrow 0$ and $T\uparrow\infty$ everywhere and use
  $\check{c}_{h,\lambda},$ $\check\varkappa,$ Proposition~\ref{bw_av} instead of $\hat{c}_{h,T}$,$\hat{\varkappa}$,
  Proposition~\ref{av_bw}.
  \bo

  Note that the conditions of the corollary regarding the precision function  can be relaxed.
  For example, in the case where  the limit $S_*$ of $\avm_T$ (or $\avp_T$)  as $T\uparrow \infty$ exists and is uniform   on $\Omega$,
   in Corollary~\ref{maintheoremgeneral1}, it is sufficient to provide the existence of
    a precision  function for $\av_T(T\uparrow \infty)$
    and common precisions for sets
     $\{\hat{c}^{\avm}_{h,T}(T\uparrow \infty)\,|\,h>0\}$,
         $\{\hat{c}^{\avp}_{h,T}(T\uparrow \infty)\,|\,h>0\}.$
    Moreover,  common precisions for sets
     $\{\hat{c}^{\avm}_{h,T}(T\uparrow \infty)\,|\,h>0\}$,
         $\{\hat{c}^{\avp}_{h,T}(T\uparrow \infty)\,|\,h>0\},$
     $\{\hat{c}^{S_*}_{h,T}(T\uparrow \infty)\,|\,h>0\}$
     can only exist simultaneously, i.e., it is sufficient to check any one of them.

{\bf On Tauberian theorem for abstract control system.}

We can obtain
the Tauberian theorem for abstract control system \cite{barton}.
Following \cite{barton}, assume the sets $\Omega$, $\mm{K}$ to be given; moreover, assume $\mm{K}$ to be closed under concatenation and $\Gamma(\omega)$, defined as above, to be non-empty for all $\omega\in\Omega$. Set $\mB\rav\{\mm{K}\}.$
Consider all the possible selectors  of the multivalued
mapping $\omega\mapsto\Gamma(\omega)$; let
$\mA$ be the set of all
possible  images  of  these  mappings.
The closedness of $\mm{K}$ under concatenation implies
the same for $\mA$,$\mB$.
We can directly check that, for  a bounded payoff, the
corresponding games have  saddle points and $\epsi$-optimal players' strategies for a positive
$\epsi.$ Therefore, by Corollary~\ref{maincorollary}, we get
the Tauberian theorem for abstract control system \cite{barton}.

{\bf On Tauberian theorem for normal form games. Proof of Theorem~\ref{normal1}}

 For all $l\in\mathfrak{L}$, $m\in\mathfrak{M},$ define $A_l\rav
 \big\{z[\omega,l,m']\,\big|\,\omega\in\Omega,m'\in\mathfrak{M}\big\},\quad
 B_m\rav
 \big\{z[\omega,l',m]\,\big|\,\omega\in\Omega,l'\in\mathfrak{L}\big\}.$
 Set $\mA\rav\{A_l\,|\,l\in\mathfrak{L}\}$,$\mB\rav\{B_m\,|\,m\in\mathfrak{M}\}.$
 It is easy see that $A_l\cap B_m\cap \Gamma(\omega)=\{z[\omega,l,m]\}\neq\pust$ for all $l\in\mathfrak{L}$, $m\in\mathfrak{M},$ $\omega\in\Omega.$
 Thus,  condition $(\ct{P})$  holds for $\mA,\mB.$

 For each mapping $\zeta:\Omega\to\mA$, for each $\omega\in\Omega$, there exists  a  $l^\omega\in\ct{L}(\omega)$ such that
 $\zeta(\omega)\cap \Gamma(\omega)=A_{l^\omega}\cap \Gamma(\omega).$ Define $l^*\in\mathfrak{L}$ by the rule $l^*(\omega)\rav l^\omega(\omega)$ for all $\omega\in\Omega.$ Now, $A_{l^*}\in\mA$ satisfies
 $\zeta(\omega)\cap \Gamma(\omega)=A_{l^*}\cap \Gamma(\omega)$ for all $\omega\in\Omega.$
 So, condition $(\omega)$  holds for $\mA.$ Analogously,  it is easy to prove that $\mB$ is the same. Thus, the result of Lemma~\ref{firstlemm1} holds.

 Moreover, comparing \rref{newa} and \rref{310}, we obtain that,
  for all $\tau>0,$ for all $A_{l'},A_{l''}\in\mA$,
 there exists $A_l\in\mA$ such that $A_l\cap\Gamma(\omega)=(A_{l'}\diamond_\tau A_{l''})\cap\Gamma(\omega)$ for all  $\omega\in\Omega$. So, $\mA$ is closed with respect to concatenation.
 Analogously,  \rref{newb} and \rref{310} imply  that $\mB$ is the same.

 Now, comparing relations \rref{2900a},\rref{2900b} with relations \rref{280},\rref{580},\rref{248}, we obtain
 $\avm\equiv\avg\equiv\avp,\bwm\equiv\bwg\equiv\bwp,$ i.e., for all $\lambda,T>0,$
 for each of the following payoffs
             $\av_T$,
             $\bw_\lambda$,
             the corresponding  games have saddle points.
 Comparing relations \rref{2950a},\rref{2950b} with relations \rref{533},\rref{633}, we obtain
 the same for payoffs $\hat{c}^{\avm}_{h,T}$, $\check{c}^{\bwm}_{h,\lambda}$ for each $\lambda,T,h>0.$

 By Lemma~\ref{firstlemm1},
 for each $\lambda,T,h>0$, for each of the following payoffs
             $\av_T$, $\hat{c}^{\avm}_{h,T},$
             $\bw_\lambda$, $\check{c}^{\bwm}_{h,\lambda}$,
             the corresponding  games have
$\epsi$-optimal player's strategies from $\mA,\mB$ respectively for all $\epsi>0$.

All conditions of Theorem~\ref{maintheoremgeneral} are verified. By this theorem we obtain what was needed. \bo

{\bf Proof of Theorem~\ref{maintheoremgeneralnew}.}

Note that condition~$(\ct{C})$ implies the result of Lemma~\ref{firstlemm2}.
Therefore, $\cym[c]\leq\cyp[c]$ for every bounded payoffs $c.$
Then,  it is sufficient to prove \rref{889}  for two  payoff families:
$\av_T(T\uparrow\infty)$ and $\bw_\lambda(\lambda\downarrow 0)$.

The condition $(\omega)$ implies the result of  Lemma~\ref{firstlemm1}. By Remark~\ref{777},
two payoff families  $\av_T(T\uparrow \infty)$,
             $\bw_\lambda (\lambda\downarrow 0)$
              have  precisions.

By condition, one of both limits,
either the common limit of $\avm_T$ and $\avp_T$ as $T\uparrow\infty,$
or the common limit of $\bwm_\lambda$ and $\bwp_\lambda$ as $\lambda\downarrow 0,$
 exists and
        is uniform  in $\omega\in\Omega.$ Denote this limit by $S_*.$
The corresponding   payoff family has a precision. Then,
by Remark~\ref{888}, $S_*$ is a protected asymptotic guarantee for  lower and upper game families with this payoff family.

Now,  condition~$(\tau)$ implies the result of  Lemma~\ref{crucial} for this payoff family.
By Remark~\ref{888},  $S_*$
is a  subsolution for the corresponding  game family.
Then, by Corollary~\ref{maincorollary}, $S_*$ is
 a common protected asymptotic guarantee  for both
       lower game families with payoffs $\av_T(T\uparrow\infty)$ and with payoffs $\bw_\lambda(\lambda\downarrow 0)$.

Because $\mB$ satisfies the same set of conditions as $\mA$ does, the assumptions of this
 Theorem with respect to lower and upper games are totally symmetric.
Then,  this limit $S_*$ is
 a common protected asymptotic guarantee  for
       upper game families with payoffs $\av_T(T\uparrow\infty)$ and payoffs $\bw_\lambda(\lambda\downarrow 0)$.
This implies that \rref{889} holds.
\bo

As follows from this proof, condition $(\omega)$ can be relaxed:
    \begin{corollary}\label{999}
              Let $\mA$ and $\mB$ satisfy conditions $(\ct{C})$,$(\diamond)$,$(\tau)$.

 If payoff family $\av_T(T\uparrow\infty)$ has a precision and the  limits of $\avm_T,\avp_T$ (in \rref{621}) exist as $T\uparrow\infty$,
 are uniform on $\Omega,$ and coincide, then
all limits in \rref{621} exist, are uniform on $\Omega,$ and coincide.

If payoff family $\bw_\lambda(\lambda\downarrow 0)$ has a precision and the  limits of $\bwm_\lambda,\bwp_\lambda$ (in \rref{621}) exist as $\lambda\downarrow 0$,
 are uniform on $\Omega,$ and coincide, then
all limits in \rref{621} exist, are uniform on $\Omega,$ and coincide.

Moreover, in each of this cases,
  $S_*$ is   a common protected asymptotic guarantee  for
     lower and  upper game families with payoffs $\av_T(T\uparrow 0)$ and $\bw_\lambda(\lambda\downarrow 0)$.
   \end{corollary}
\begin{remark}\label{555}
   Under conditions of Corollary~\ref{999}, we can apply Remark~\ref{111} from Appendix~\ref{AA}. Then, by
    Remark~\ref{222} from Appendix~\ref{A} and Remark~\ref{333} from Appendix~\ref{B}, the strategies that protect a guarantee $S_*$
   for one of the game families (with payoffs $\av_T(T\uparrow 0)$ or $\bw_\lambda(\lambda\downarrow 0)$),
   are expressed (see  \rref{UT},\rref{Ul}) through
    such strategies for the other game family.
\end{remark}

{\bf On Tauberian theorem for normal form games. Proof of Theorem~\ref{normal2}.}

 For all $l\in\mathfrak{L}$, $m\in\mathfrak{M},$ define $A_l\rav
 \big\{z[\omega,l,m']\,\big|\,\omega\in\Omega,m'\in\mathfrak{M}\big\},\quad
 B_m\rav
 \big\{z[\omega,l',m]\,\big|\,\omega\in\Omega,l'\in\mathfrak{L}\big\}.$
 Set $\mA\rav\{A_l\,|\,l\in\mathfrak{L}\}$,$\mB\rav\{B_m\,|\,m\in\mathfrak{M}\}.$
 It is easy see that $A_l\cap B_m\cap \Gamma(\omega)=\{z[\omega,l,m]\}\neq\pust$ for all $l\in\mathfrak{L}$, $m\in\mathfrak{M},$ $\omega\in\Omega.$
 Thus,  conditions $(\ct{P})$ and $(\ct{C})$ hold for $\mA,\mB.$

 For each mapping $\zeta:\Omega\to\mA$, for each $\omega\in\Omega$, there exists  a  $l^\omega\in\ct{L}(\omega)$ such that
 $\zeta(\omega)\cap \Gamma(\omega)=A_{l^\omega}\cap \Gamma(\omega).$ Define $l^*(\omega)\rav l^\omega(\omega).$ Now,
 $\zeta(\omega)\cap \Gamma(\omega)=A_{l^*}\cap \Gamma(\omega)$ for all $\omega\in\Omega.$ Thus,
  condition $(\omega)$ holds for $\mA.$  Analogously,  it is easy prove that $\mB$ is the same.

 For all  $\tau>0$, for all $l',l''\in\mathfrak{L}$, $\omega\in\Omega$,
\begin{eqnarray*}
A_{l'}\diamond_\tau A_{l''}&\ravref{310}&\big\{z'\diamond_\tau z''\,\big|\,z'\in A_{l'},z''\in A_{l''},z'(\tau)=z''(0)\big\}\\
&=&\big\{z[\omega,l',m']\diamond_\tau z[\omega',l'',m'']\,|\,m',m''\in\mathfrak{M}, \omega'=z[\omega,l',m'](\tau)\big\}\\
&\ravref{buka}&\big\{z[\omega,l'\diamond_\tau l'',m'\diamond_\tau m'']\,|\,m',m''\in\mathfrak{M}\big\}\\
&\ravref{byka}&\big\{z[\omega,l'\diamond_\tau l'',m]\,|\,m\in\mathfrak{M}\big\}=A_{l'\diamond_\tau l''}\in\fr{A}.
\end{eqnarray*}
 So,  $\mA$ is closed with respect to concatenation.
 Moreover, by \rref{byka}, for all  $l'\in\mathfrak{L}$, there exists $l''\in\mathfrak{L}$ satisfying $l'\diamond_\tau l''=l'.$
 Then, $A_{l'}\diamond_\tau A_{l''}=A_{l'}.$  So,  $\mA$ is closed under backward shift.
 Analogously,  it is easy prove that $\mB$ is the same.

 Now, comparing relations \rref{naa},\rref{nab} with relations \rref{280},\rref{580},\rref{248}, we see that
 all conditions of Corollary~\ref{999} are verified. By this corollary we obtain what was required. \bo

\section{Games for discrete setting}
\label{normalform}

{\bf Game statement.}

Recall that $\mm{N}=\{1,2,3,\dots\}.$
Let  a set $\bar{\Omega}$ be nonempty.
   Let $\bar{\mm{K}}$
   be a set of admissible   processes $\bar{y}:\mm{N}\cup\{0\}\to\Omega$.
    For all $\bar{\omega}\in\bar{\Omega},$
    consider a
   set  $\bar{\Gamma}(\bar{\omega})\rav\{\bar{y}\in\bar{\mm{K}}\,|\,\bar{y}(0)=\bar{\omega}\}$
of processes $\bar{y}\in\bar{\mm{K}}$
   with $\bar{y}(0)=\bar{\omega}$. Let these sets be also non-empty for all  $\bar{\omega}\in\bar{\Omega}.$

Again,  the first player wishes  to maximize
  a payoff function $c:\bar{\mm{K}}\to \mm{R}$; the second player wishes to minimize it.
   Players are given the sets $\mA,\mB$ of subsets of $\mm{\bar{K}}$, their sets of strategies.
To ensure that the lower and upper games are well-defined, assume that,
for each $\bar\omega\in\bar\Omega$, for all $\bar{A}\in\bar{\mA}$,$\bar{B}\in\bar{\mA}$, sets $\bar{A}\cap\bar{B}\cap \bar\Gamma(\bar\omega)$
are non-empty.

   Fix a function  $\bar{g}:\bar\Omega\to [0,1].$
 For all $\mu\in(0,1), n\in\mm{N}$, define payoffs $\bar{v}_n, \bar{w}_\mu$ as follows:
 \begin{eqnarray*}
 \bar{\av}_n(\bar{y})=\frac{1}{n}\sum_{t=0}^{n-1} \bar{g}(\bar{y}(t))\in [0,1],\quad
 \bar{\bw}_\mu(\bar{y})=\mu\sum_{t=0}^\infty (1-\mu)^{i} \bar{g}(\bar{y}(t))\in[0,1]\qquad \forall \bar{y}\in\bar{\mm{K}}.
\end{eqnarray*}

 For simplicity, assume that, both for the  payoffs $\bar{\av}_n$ and payoffs $\bw_\mu$,
 the corresponding games have their saddle points, i.e.,  the following definitions are valid:
 \begin{subequations}
  \begin{eqnarray}
\bavg_n(\bar\omega)&=&
\sup_{\bar{A}\in\bar{\mA}}\inf_{\bar{y}\in\bar{A}\cap\bar\Gamma(\bar\omega)}\bar{\av}_n(\bar{y})=
\inf_{\bar{B}\in\bar{\mB}}\sup_{\bar{y}\in\bar{B}\cap\bar\Gamma(\bar\omega)}\bar{\av}_n(\bar{y}),
\qquad \forall n\in\mm{N},\bar\omega\in\bar\Omega,
\label{1576}
\\
\bbwg_\mu(\bar\omega)&=& \sup_{\bar{A}\in\bar{\mA}}\inf_{\bar{y}\in\bar{A}\cap\bar\Gamma(\bar\omega)}\bar{\bw}_\mu(\bar{y})=
\inf_{\bar{B}\in\bar{\mB}}\sup_{\bar{y}\in\bar{B}\cap\Gamma(\bar\omega)}\bar{\bw}_\mu(\bar{y})
\qquad \forall \mu\in(0,1),\bar\omega\in\bar\Omega.\label{1577}
  \end{eqnarray}
 \end{subequations}
  We will also need a condition to guarantee the existence of near-optimal strategies for the players.
  For all $\epsi>0$, or, in the games with the payoffs $\bw_\mu$
  $(\mu\in(0,1))$ for sufficiently small $\mu$, or, respectively, the games with the   payoffs $\av_n$
  $(n\in\mm{N})$ for sufficiently large $n$,
  each player has $\epsi$-optimal strategies, i.e., for all $\epsi>0$ there exists $N\in\mm{N}$ such that
 \begin{subequations}
  \begin{eqnarray}
\forall n>N\
\exists\bar{A}^{n,\epsi}\in\mA,\bar{B}^{n,\epsi}\in\mB
\
\forall\bar\omega\in\bar\Omega\quad
\inf_{\bar{y}\in\bar{A}^{n,\epsi}\cap\bar\Gamma(\bar\omega)}\bar{\av}_n(\bar{y})+\epsi\!\geq\!
\bavg_n(\bar\omega)\!\geq\!\!\sup_{\bar{y}\in\bar{B}^{n,\epsi}\cap\bar\Gamma(\bar\omega)}\bar{\av}_n(\bar{y})-\epsi,\quad
\label{1578}\\
\forall \mu\in(0,{\textstyle\frac{1}{N}})\
\exists\bar{A}^{\mu,\epsi}\in\mA,\bar{B}^{\mu,\epsi}\in\mB
\
\forall\bar\omega\in\bar\Omega\
\inf_{\bar{y}\in\bar{A}^{\mu,\epsi}\cap\bar\Gamma(\bar\omega)}\bar{\bw}_\mu(\bar{y})+\epsi\!\geq\!
\bbwg_\mu(\bar\omega)\!\geq\!\!\sup_{\bar{y}\in\bar{B}^{\mu,\epsi}\cap\bar\Gamma(\bar\omega)}\bar{\bw}_\mu(\bar{y})-\epsi.\quad
\label{1579}
  \end{eqnarray}
 \end{subequations}

 Let, for each $n\in\mm{N}$, for all $\bar{y}',\bar{y}''\in\bar{\mm{K}}$ with $\bar{y}'(n)=\bar{y}''(0)$,
  define their concatenation $\bar{y}'\diamond_\tau \bar{y}'':\mm{N}\to \bar{\Omega}$ by the rule:
$$(\bar{y}'\diamond_n \bar{y}'') (i)\rav \left\{
 \begin{array} {rcl}        \bar{y}'(i),       &\mathstrut&     i\leq n,i\in\mm{N}\cup\{0\};\\
                            \bar{y}''(i-n),
                                         &\mathstrut&    i>n,i\in\mm{N}.
 \end{array}            \right.$$
 So,  for all subsets $\bar{A}',\bar{A}''\subset\bar{\mm{K}}$  and a time $n\in\mm{N}$, define
  their concatenation $\diamond_n$ by
\begin{eqnarray*}
  \bar{A}'\diamond_n \bar{A}''&\rav&
  \{\bar{y}'\diamond_n \bar{y}''\,|\,\bar{y}'\in\bar{A}',\bar{y}''\in \bar{A}'',\bar{y}''(n)=\bar{y}'(0)\}.
\end{eqnarray*}

 Let us say that $\bar{\mA}$   is  closed under concatenation  and backward shift  with integral time points if,
 for a natural $n\in\mm{N}$ and for  $\bar{A},\bar{A}'\in\mA$, we have $\bar{A}\diamond_n \bar{A}'\in\bar{\mA}$ and
  $\bar{A}=\bar{A}\diamond_n \bar{A}''$ for some $\bar{A}''\in\bar{\mA}.$

{\bf Tauberian theorem for games with discrete time.}

\begin{theorem}
\label{normtauberian}
 Let $\bar{\mA},\bar{\mB}$ be sets of  strategies,
 and let they  be  closed  under concatenation  and backward shift  with integral time points.
 Assume also that $\bar{A}\cap \bar{B}\cap \bar{\Gamma}(\bar{\omega})$ for all $\bar{\omega}\in\bar\Omega,\bar{A}\in\bar{\mA},\bar{B}\in\bar{\mB}.$

 Let, for all $n\in\mm{N},$ $\mu\in(0,1)$, \rref{1576}-\rref{1577} hold,
 i.e., let all  games with payoffs $\bar{\av}_n$,$\bar{\bw}_\mu$ have saddle points.

 If \rref{1578} holds and the limit of $\bavg_n$ exists as $n\to\infty$  and is
 uniform in $\bar\omega\in\bar\Omega,$ then
both limits
\begin{eqnarray} \label{v}
 \lim_{n\to\infty} \bavg_n(\bar\omega)=
  \lim_{\mu\downarrow 0} \bbwg_\mu(\bar\omega)
  \qquad\forall\bar\omega\in\bar\Omega
\end{eqnarray}
 exist, are uniform on $\bar\Omega,$ and coincide; moreover, \rref{1578},\rref{1579} hold.

 On the other hand, if \rref{1579} holds and  the limit of $\bbwg_\mu$ exists as $\mu\downarrow 0$, and is
 uniform in $\bar\omega\in\bar\Omega,$ then both limits in \rref{v}
 exist, are uniform on $\bar\Omega$, and coincide; moreover, \rref{1578},\rref{1579} hold.
\end{theorem}

{\bf Proof of Theorem~\ref{normtauberian}. Reduction to general statement.}

 To each $t\geq 0$, assign the integral number $\lfloor t\rfloor$, the greatest number (from $\mm{N}\cup\{0\}$) not surpassing  $t.$ Now, for all $t\geq0$, $\ldr t\rdr \rav t-\lfloor t\rfloor\in[0,1)$ is the fractional part of $t$.

  Define $S\rav[0,1)$,$\Omega\rav \bar{\Omega}\times S=\bar{\Omega}\times [0,1)$.
  To each $\bar{y}\in\bar{\mm{K}},s\in S$, assign $z[\bar{y},s]:\mm{T}\to\Omega$
  by the rule: $z[\bar{y},s](t)=(\bar{y}(\lfloor t+s\rfloor),\ldr t+s\rdr)$ for all $t\in\mm{T}.$
   Set
    $$\mm{K}\rav \{z[\bar{y},s]\,|\,\bar{y}\in\bar{\mm{K}},s\in S\}.$$

  Define $g(\bar{\omega},s)\rav\bar{g}(\bar{\omega})$  for all $\omega=(\bar{\omega},s)\in\Omega$.
  Then, for all $z[\bar{y},s]\in\mm{K}$, $g\big(z[\bar{y},s](t)\big)=\bar{g}\big(\bar{y}(\lfloor t+s\rfloor)\big).$

  Denote by $\mm{A}$ and $\mm{B}$, respectively, the sets of all possible mappings   of  $S\ni s\to\xi(s)\in\bar\mA$ and
   $S\ni s\to\xi(s)\in\bar\mB.$
    For all $\xi\in\mm{A}\cup\mm{B}$, designate
   \begin{eqnarray*}
   Z[\xi]\rav\{z[\bar{y},s]\,|\, s\in{S},\bar{y}\in \xi(s)\},\quad
    \mA\rav\{Z[\xi]\,|\,\xi\in\mm{A}\},\quad
    \mB\rav\{Z[\xi]\,|\,\xi\in\mm{B}\}.
   \end{eqnarray*}

  Now, for all $(\bar\omega,s)\in{\Omega},A=Z[\xi_A]\in\mA,
  B=Z[\xi_B]\in\mB$, we obtain
\begin{eqnarray*}
    \Gamma(\bar\omega,s)&=&\{z[\bar{y},s']\,|\,\exists \bar{y}\in\bar{K},s'\in S,\bar{y}(\lfloor s'\rfloor)=\bar{\omega},\ldr s'\rdr=s\}=
    \{z[\bar{y},s]\,|\,\bar{y}\in \bar\Gamma(\bar\omega)\}\\
    A\cap B\cap \Gamma(\bar\omega,s)&=&\{z[\bar{y},s']\,|\,\bar{y}\in \xi_A(s)\cap\xi_B(s)\cap\bar{\Gamma}(\bar\omega)\}\neq\varnothing.
\end{eqnarray*}
   In particular, conditions $(\ct{P})$ and $(\ct{C})$ hold for $\mA,\mB$.

{\bf Proof of Theorem~\ref{normtauberian}. Closedness with respect to concatenation and backward shift.}

  It is sufficient to examine the set $\bar{\mA}.$

  To simplify the notation,
  for all $\bar{y}',\bar{y}''\in\bar{\mm{K}}$ with $\bar{y}'(0)=\bar{y}''(0)$, set
  $\bar{y}'\diamond_0\bar{y}''\rav\bar{y}'';$
  now,
    $\bar{A}'\diamond_0 \bar{A}''= \bar{A}''$
for all $\bar{A}',\bar{A}''\in\bar{\mA}.$

  Remember that, for all $\tau>0$, one has $z'=z[\bar{y}',s'],z''=z[\bar{y}'',s'']\in\mm{K}$;
   their concatenation $z'\diamond_\tau z''$ is defined  if   $z'(\tau)=z''(0).$
  Now, $z'(\tau)=z''(0)$ iff
  $\ldr\tau+s'\rdr=\ldr s''\rdr,\ \bar{y}'(\lfloor \tau+s'\rfloor)=\bar{y}''(\lfloor s''\rfloor)$
   hold;
   by $s''\in S$, it is now equivalent to the pair of equalities:
   $\ldr\tau+s'\rdr=s'',\bar{y}'(\lfloor \tau+s'\rfloor)=\bar{y}''(0).$
   Also in this case,  $\ldr\tau+s'+t\rdr=\ldr s''+t\rdr$,$\lfloor\tau+s'+t\rfloor=\lfloor s''+t\rfloor$ for all $t\geq 0.$
   Then,
   we see that
   $z[\bar{y}',s']\diamond_\tau z[\bar{y}'',s'']$ is well-defined iff
   $z[\bar{y}',s']\diamond_\tau z[\bar{y}'',s'']=(\bar{y}'\diamond_{\lfloor\tau+s'\rfloor}\bar{y}'',s'),$
      $s''=\ldr\tau+s'\rdr$ and $\bar{y}'(\lfloor \tau+s'\rfloor)=\bar{y}''(0)$ hold, i.e.,
   $\bar{y}'\diamond_{\lfloor \tau+s'\rfloor}\bar{y}''$ is well-defined, and $s''=\ldr\tau+s'\rdr$.

  So,
 for all
 $A'=Z[\xi'],A''=Z[\xi'']\in\mA$, $\tau>0,$ define
  $\xi\in\mm{A}$ by the rule  $\xi(s)\rav\xi'(s)\diamond_{\lfloor\tau+s\rfloor}\xi''(\ldr\tau+s\rdr)$
 for all $s\in{S}.$ Then, we obtain
 \begin{eqnarray}
 A'\diamond_\tau A''&=&\big\{z'\diamond_\tau z''\,\big|\,z'\in A',z''\in A'',z'(\tau)=z''(0)\,\big\}\nonumber\\
 &=&\big\{z[\bar{y}',s']\diamond_\tau z[\bar{y}'',s'']\,\big|\,s'\in S,\bar{y}'\in\xi'(s'),s''=\ldr\tau+s'\rdr,\bar{y}''\in \xi''(s''),\bar{y}'(\lfloor \tau+s'\rfloor)=\bar{y}''(0)\big\}\nonumber\\
 &=&\big\{z[\bar{y}'\diamond_{\lfloor\tau+s'\rfloor} \bar{y}'',s']\,\big|\,
 s'\in S,\bar{y}'\in\xi'(s'),\bar{y}''\in \xi''(\ldr\tau+s'\rdr),\bar{y}'(\lfloor \tau+s'\rfloor)=\bar{y}''(0)\big\}\nonumber\\
 &=&\big\{z[\bar{y},s]\,\big|\,s\in S,\bar{y}\in \xi'(s)\diamond_{\lfloor\tau+s\rfloor}\xi''(\ldr\tau+s\rdr)\big\}\label{1430}\\
 &=&\big\{z[\bar{y},s]\,\big|\,s\in S,\bar{y}\in \xi(s)\big\}
 =Z[\xi]\in\mA,\nonumber
 \end{eqnarray}
 Thus, condition $(\diamond)$ proved.

 Consider  arbitrary
 $A=Z[\xi]$, $n\in\mm{N}$, $\tau\in(n-1,n]$. For all $s\in{S}$, by condition,
 $\xi(s)=\xi(s)\diamond_{\lfloor\tau+s\rfloor} \bar{A}^{(s)}$ for some $\bar{A}^{(s)}\in\bar{\mA}$.
 Define $\xi^*\in\mm{A}$ by the rule
 $\xi^*(\ldr\tau+s\rdr)=\bar{A}^{(s)}$ for all $s\in S$ (i.e.
 $\xi^*(s')=\bar{A}^{(s'+n-\tau)}$ for $s'\in [0,\tau-n+1)$,
 $\xi^*(s')=\bar{A}^{(s'+n-1-\tau)}$ for $s'\in [\tau-n+1,1)$). Now, $\xi(s)=\xi(s)\diamond_{\lfloor\tau+s\rfloor} \xi^*(\ldr\tau+s\rdr)$ for all $s\in{S}.$
 Then, $A^*\rav Z[\xi^*]\in\mA$ satisfies
 \begin{eqnarray*}
 A\diamond_\tau A^*&\ravref{1430}&\big\{z[\bar{y},s]\,\big|\,s\in S,\bar{y}\in \xi(s)\diamond_{\lfloor\tau+s\rfloor}\xi^*(\ldr\tau+s\rdr)\big\}\\
 &=&\big\{z[\bar{y},s]\,\big|\,s\in S,\bar{y}\in \xi(s)\big\}=Z[\xi]=A.
 \end{eqnarray*}
 Thus, condition $(\tau)$ is also proved.


{\bf Proof of Theorem~\ref{normtauberian}. Payoffs' comparison.}

  Consider the following function $\mu:\mm{R}_{>0}\to (0,1)$:
  $\mu(\lambda)=1-e^{-\lambda}$  for all $\lambda>0$.
  Note that $\mu(0+)=0+.$

  Fix a process  $z\in\mm{K}$. Now, $z=z[\bar{y},s]$ for some $\bar{y}\in\bar{\mm{K}}, s\in S.$
  In addition, for all $\bar{y}\in\bar{\mm{K}}, s\in S$ we can find such $z.$

  Define $z'\rav z[\bar{y},0]$.
  Then, $z\equiv z'_s$ (see \rref{296}) and
  $g(z(t-s))=\bar{g}(\bar{y}(\lfloor t\rfloor))=g(z'(t))$
   for all $t\geq s.$

  \begin{subequations}
    In addition,  for all $n\in\mm{N}$,
\begin{eqnarray}
    \bar{\av}_n(\bar{y})=\frac{1}{n}\sum_{t=0}^{n-1}
     \bar{g}(\bar{y}(t))=
       \frac{1}{n}\sum_{i=0}^{n-1}\int_{i}^{i+1}g(z'(r))\,dr=\frac{1}{n}\int_{0}^{n}g(z'(r))\,dr={\av}_{n}(z').
       \label{1730}
\end{eqnarray}
  Moreover,
  $\int_{t}^{t+1}\lambda e^{-\lambda r}\,dr=e^{-\lambda t}-e^{-\lambda (t+1)}=(1-\mu(\lambda))^{t}\mu(\lambda)$
  for all $t\geq 0,\lambda>0$; now,
\begin{eqnarray}
\label{1740}
 \bar{\bw}_{\mu(\lambda)}(\bar{y})&=&\sum_{i=0}^\infty \mu(\lambda)(1-\mu(\lambda))^{i}\int_{i}^{i+1} g(z'(r))\,dr=   \int_{0}^{\infty} \lambda e^{-\lambda t}g(z'(r))\,dr=\bw_\lambda(z').
\end{eqnarray}
 \end{subequations}

 Define
 $\hat{\eta}(T)\rav \frac{4}{T},\quad \check{\eta}(\lambda)\rav 2\lambda$ for all $T,\lambda>0.$
In Appendix~\ref{AA}, we prove inequalities  \rref{2228};   now,  for all $T,\lambda>0$, we obtain
 \begin{subequations}
\begin{eqnarray}
\label{important}
|\bar{\av}_{\lfloor T \rfloor+1}(\bar{y})-\av_T(z)|\ravref{1730}|\av_{\lfloor T \rfloor+1} (z')-\av_T(z'_s)|\leqref{2228}
 \hat{\eta}(T)&\to& 0\quad\textrm{ as }T\uparrow \infty;\\
|\bar{\bw}_{\mu(\lambda)}(\bar{y})-\bw_\lambda(z)|\ravref{1740}|\bw_{\lambda} (z')-\bw_{\lambda}(z'_s)|\leqref{2228} \check{\eta}(\lambda)&\to& 0\quad\textrm{ as }\lambda\downarrow 0.\label{important1}
\end{eqnarray}
 \end{subequations}

 Consider the case where \rref{1578} holds and the
 limit of $\bavg_{n}$ (in \rref{v}) exists and is uniform on $\bar\Omega$.
   Then, \rref{important} implies that limits of $\avm_T$ and $\avp_T$ (from \rref{621}) exist and are uniform on $\Omega$. Moreover, by \rref{important}, \rref{1578}  implies that
  the  payoff family $\av_T(T\uparrow\infty)$ has the precision $2\hat{\eta}.$
  Now, by Corollary~\ref{999}, all limits in~\rref{621} exist, are uniform on $\Omega,$ and coincide.
  Then, by \rref{important}-\rref{important1},
all limits in \rref{v}
 exist, are uniform on $\bar\Omega,$ and coincide. Moreover, $S_*$ is a common protected asymptotic guarantee  for
     lower and  upper game families with payoffs $\bw_\lambda(\lambda\downarrow 0)$.
By \rref{important1} and $\check{\eta}(\lambda)\to 0+$, it follows \rref{1579}.

 The case where
 the limit of $\bbwg_{\mu}$  exists and is uniform on $\bar\Omega$
 is analyzed analogously.
  \bo




\setcounter{section}{0}%
\def\thesection{\Alph{section}}%
\section{Auxiliary facts. The proofs of lemmas.}
\label{AA}

 \begin{subequations}
       Let $\nu_\gamma(\gamma\to\gamma_*)$ be one of the
 payoff families, either $\av_T(T\uparrow\infty),$ or
       $\bw_\lambda(\lambda\downarrow 0)$,

For all $h>0,$ under $\nu_\gamma=\av_T$, set
  $\gamma=T$,$\gamma_h=T+h$,  $\sigma_{h,\gamma}=\frac{T}{T+h}=\frac{\gamma}{\gamma_h},$
  $\rho_\gamma|_{[0,\gamma]}\equiv\frac{1}{\gamma},$
  $\rho_\gamma|_{(\gamma,\infty)}\equiv 0;$
 under $\nu_\gamma=\bw_\lambda$, set
  $\gamma=\lambda$,$\gamma_h=\lambda$, $\sigma_{h,\gamma}=e^{-\lambda h},$ $\rho_\gamma(t)=\lambda e^{-\lambda t}$ for all $t\geq 0.$

For a function $U:\mm{R}_{>0}\times\Omega\mapsto\mm{R}$, for all $h>0$, we can define a payoff family by the following rule:
  \begin{eqnarray}
 \label{20002000}
c^{U}_{h,\gamma}(z)\rav 
\int_0^h \rho_{\gamma_h}(t)g(z(t))\,dt+\sigma_{h,\gamma} U_{\gamma}(z(h))\quad\forall z\in\mm{K}\quad(\gamma\to\gamma_*).
  \end{eqnarray}
Then, by \rref{533},\rref{633}, $c^{U}_{h,T}\equiv\hat{c}^{U}_{h,T}$ if $\nu_\gamma\equiv \av_T$,
 and ${c}^{U}_{h,\gamma}\equiv \check{c}^{U}_{h,\lambda}$ if $\nu_\gamma\equiv \bw_\lambda$.

 Note several useful properties:
 for all $\gamma,h>0,T\geq 0,z\in\mm{K},z'\in\Gamma(z(h))$
  \begin{eqnarray}
  \rho_{\gamma_h}(h+T)=\sigma_{h,\gamma}\rho_{\gamma}(T),  \quad
   \nu_\gamma(z)=\int_0^\infty \rho_{\gamma}(t) g(z(t))\,dt;
 \label{2000g}\\
 \label{2000m}
  \max(0,1-h\rho_\gamma(0))<\sigma_{h,\gamma}<1,\quad\int_{0}^\infty \rho_{\gamma}(t)\ dt=1;\\
  \label{2000b}
  h\rho_\gamma(0)\geq\int_{0}^h \rho_{\gamma_h}(t)z(t)\ dt\ravref{2000g}
  \nu_{\gamma_h}(z)-\int_{0}^\infty \sigma_{h,\gamma}\rho_{\gamma}(t)z(t+h)\ dt\ravref{296}
  \nu_{\gamma_h}(z)-\sigma_{h,\gamma}\nu_{\gamma}(z_h);\\
\label{2000h}
  c^{U}_{h,\gamma}(z)\ravref{20002000}
  \int_0^h \rho_{\gamma_h}(t)g(z(t))\,dt+\sigma_{h,\gamma} U_{\gamma}(z(h))\ravref{20002000}c^U_{h,\gamma}(z\diamond_h z'). \end{eqnarray}

We claim also that
 for every
 $z\in\mm{K},s,r\in[0,1],\lambda,T> 0$, one has
\begin{eqnarray}
\label{2228}
 |\bw_{\lambda}(z)-\bw_{\lambda}(z_s)|\leq 2\lambda,\qquad
 |\av_{T+r}(z)-\av_{T}(z_s)|\leq\frac{4}{T}.
\end{eqnarray}
 \end{subequations}
 \begin{subequations}
Indeed,
       let $\nu_\gamma(\gamma\to\gamma_*)$ be one of the
 payoff families, either $\av_T(T\uparrow\infty),$ or
       $\bw_\lambda(\lambda\downarrow 0)$.
 For all $z\in\mm{K},\gamma>0,h\in[0,1)$, we obtain
\begin{eqnarray}
\label{2000w}
|\nu_{\gamma_h}(z)-\nu_{\gamma}(z_h)|\!\ravref{2000b}\!\Big|\int_{0}^h \rho_{\gamma_h}(t)z(t)\, dt\Big|+(1-\sigma_{h,\gamma})\nu_{\gamma}(z_h)\!\leqref{2000b}\!
 h\rho_\gamma(0)+|1-\sigma_{h,\gamma}|\!\leqref{2000m}\!2\varrho_\gamma(0).
\end{eqnarray}
Now, to get the first inequality from \rref{2228}, it is sufficient to make a substitution $\gamma=\lambda,\nu_\gamma=\bw_\lambda,h=r$ in \rref{2000w}.

 It remains to consider the case $\gamma=T,\nu_\gamma=\av_T.$
 Recall that, in this case, $\rho_\gamma|_{[0,\gamma]}\equiv\frac{1}{\gamma},$
  $\rho_\gamma|_{(\gamma,\infty)}\equiv 0.$ Now, for all $h\in[0,1]$, we have $\gamma_h=\gamma+h,$
\begin{eqnarray}
  \label{2000mm}
  \int_{0}^\infty |\rho_{\gamma}(t)-\rho_{\gamma+h}(t)|\ dt\leq
    \Big(\frac{1}{\gamma}-\frac{1}{\gamma+h}\Big)\int_{0}^\gamma\,dt+
    \frac{1}{\gamma+h}\int_{\gamma}^{\gamma+h}\,dt=\frac{2h}{\gamma+h}\leq \frac{2}{\gamma}=2\rho_{\gamma}(0).
\end{eqnarray}
Now, for all $r,s\in[0,1]$, we have $h=|\gamma_r-\gamma_s|\leq 1$. Then,
$$|\nu_{\gamma_r} (z)-\nu_{\gamma}(z_s)|
\leqref{2000w} |\nu_{\gamma_r}(z)-\nu_{\gamma_s}(z)|+2\varrho_\gamma(0)\leqref{2000mm}
4\varrho_\gamma(0),
$$
and, to get the second inequality from \rref{2228}, it suffices to substitute  $\gamma=T,\nu_\gamma=\av_T,h=r$ into the relation obtained.
So, \rref{2228} holds.
 \end{subequations}

{\bf Proof of Lemma~\ref{lastlemma}.}
            Let the set $\mA$ satisfy conditions $(\ct{P})$,$(\diamond)$.
          Let $\nu_\gamma(\gamma\to\gamma_*)$ be one of the
 payoff families, either $\av_T(T\uparrow\infty),$ or
       $\bw_\lambda(\lambda\downarrow 0)$. Fix  $\gamma,\epsi>0.$

Let $U_\gamma-\epsi$ be a protected guarantee of the lower game with the payoff $\nu_\gamma$, i.e, let
  there exist a strategy  $A^\gamma\in\mA$ satisfying
  \begin{eqnarray}
  \label{20002}
  U_\gamma(z(0))\leq \nu_{\gamma}(z)+\epsi
  \qquad\forall z\in A^\gamma.
  \end{eqnarray}
  We must prove that
  $\cym\big[c^{U}_{h,\gamma}\big]-\sigma_{h,\gamma}\epsi$ is a guarantee of the lower game with the payoff $\nu_{\gamma_h}$ for all $h>0$.

   Fix $h,\delta>0,\omega\in\Omega.$
   There exists a strategy $A'_\omega\in\mA$ satisfying
   $$\cym[c^U_{h,\gamma}](\omega)-\delta\leq
   c^U_{h,\gamma}(z)\quad \forall z\in A'_\omega\cap \Gamma(\omega).$$
   By \rref{2000h} and \rref{310}, 
   every $z\in (A'_\omega\diamond_h A^\gamma)\cap\Gamma(\omega)$
   also satisfies this inequality and $z_h\in A^\gamma$; now,    for all $z\in (A'_\omega\diamond_h A^\gamma)\cap\Gamma(\omega)$,
          \begin{eqnarray*}
       \cym[c^U_{h,\gamma}](\omega)-\delta\leq
   c^U_{h,\gamma}(z) &\ravref{20002000}&\int_0^h \rho_{\gamma_h}(t)g(z(t))\,dt+\sigma_{h,\gamma} U_{\gamma}(z(h))\\
   &\ravref{2000b}&\nu_{\gamma_h}(z)-\sigma_{h,\gamma}\nu_{\gamma}(z_h)+\sigma_{h,\gamma} U_{\gamma}(z_h(0))\\
   &\leqref{20002}&
   \nu_{\gamma_h}(z)+\sigma_{h,\gamma}\epsi.
   \end{eqnarray*}
   By condition $(\diamond)$, we have  $A'_\omega\diamond_h A^\gamma\in\mA$.
 Then, since the choice of $z\in (A'_\omega\diamond_h A^\gamma)\cap\Gamma(\omega)$ was arbitrary, we have
         \begin{eqnarray*}
       \cym[c^U_{h,\gamma}](\omega)-\sigma_{h,\gamma} \epsi-\delta   &\leq&
 \inf_{z\in (A'_\omega\diamond_h A^\gamma)\cap\Gamma(\omega)} \nu_{\gamma_h}(z)\\
   &\leq&
       \sup_{A\in\mA}\inf_{z\in A\cap\Gamma(\omega)}\nu_{\gamma_h}(z)=
   \cym[\nu_{\gamma_h}](\omega)\qquad \forall \omega\in\Omega.
   \end{eqnarray*}
  Thus,  since the positive $\delta$ was chosen arbitrarily, we have $\cym\big[c^U_{h,\gamma}]-\sigma_{h,\gamma}\epsi\leq \cym[\nu_{\gamma_h}]$ for all $h>0.$\bo

{\bf Proof of Lemma~\ref{crucial}.}

Let the set $\mA$
       satisfy conditions $(\ct{P}),(\tau)$.
       Let $\nu_\gamma(\gamma\to\gamma_*)$ be one of the
 payoff families, either $\av_T(T\uparrow\infty),$ or
       $\bw_\lambda(\lambda\downarrow 0)$. Set
       $U_\gamma(\omega)\rav\cym\big[\nu_{\gamma}](\omega)$ for all $\omega\in\Omega$,
   $\gamma>0.$ Let $U$ be a protected asymptotic guarantee for the lower game family
   with payoffs $\nu_\gamma(\gamma\to\gamma_*)$.
   We must prove that $U$ is a subsolution for this  lower game family.

   By condition, there exists a monotonic function $\varkappa:\mm{R}_{>0}\to\mm{R}_{>0}$ with
   $\varkappa(\gamma)\to 0+$ as $\gamma\to\gamma_*$ such that
   $U_{\gamma}-\varkappa(\gamma)$ is a protected guarantee for the lower game with $\nu_\gamma$ for all $\gamma>0.$

 Consider arbitrary $\gamma,h>0,A'\in\mA,z'\in A'.$
 Note that, by $(\tau)$, $A'$ may be expressed as $A'=A'\diamond_{h}A''$ for some $A''\in\mA$.
 Now,  there exists
 a ${z}''\in A''\cap \Gamma\big({z}'(h)\big)$ such that
\begin{eqnarray}
\nu_\gamma({z}'')-\varkappa(\gamma)\leq\inf_{z\in A''\cap\Gamma({z}'(h))}\nu_{\gamma}(z)\leq
\cym\big[\nu_\gamma\big]({z}'(h))=U_{\gamma}({z}'(h))
\label{2692}
\end{eqnarray}
Set
 ${z}\rav{z}'\diamond_{h}{z}''$; by ${z}_h={z}'_h={z}''$ (see \rref{296}),  we get
\begin{eqnarray*}
c^U_{h,\gamma}(z')&\ravref{2000h}&
\int_{0}^{h}\varrho_{\gamma_h}g({z}(t))\,dt+
  \sigma_{h,\gamma}U_{\gamma}({z}(h))\\
  &\geqref{2692}&
  \int_{0}^{h}\varrho_{\gamma_h}g({z}(t))\,dt+
  \sigma_{h,\gamma}\nu_\gamma({z}_h)-\sigma_{h,\gamma}\varkappa(\gamma)\\
  &\ravref{2000b}&
  \nu_{\gamma_h}({z})-\sigma_{h,\gamma}\varkappa(\gamma)\geqref{2000m}\nu_{\gamma_h}({z})-\varkappa(\gamma).
\end{eqnarray*}
So, to each $z'\in A'$, assign
${z}={z}'\diamond_{h}{z}''\in A'\diamond_{h}A''=A'$ with
$z'(0)=z(0)$ and
$c^U_{h,\gamma}(z')\geq \nu_{\gamma_h}({z})-\varkappa(\gamma).$
 Thus, we have
 $$
 \inf_{z'\in A'\cap \Gamma(\omega)}c^U_{h,\gamma}(z')+\varkappa(\gamma)
  \geq
  \inf_{z'\in A'\cap \Gamma(\omega)}\nu_{\gamma_h}(z)\qquad \forall\omega\in\Omega,A'\in\mA.$$
  Recall that,  for all $\gamma,h>0$, $U_{\gamma_h}-\varkappa(\gamma_h)$ is a protected guarantee for the lower game with payoff~$\nu_{\gamma_h}.$ 
  Then,  $U_{\gamma_h}-\varkappa(\gamma)-\varkappa(\gamma_h)$ is a protected guarantee for lower game with $c^U_{h,\gamma}$
  for all $h,\gamma>0.$ By $\varkappa(\gamma_h)\leq\varkappa(\gamma)$, we obtain that
  $U_{\gamma_h}-2\varkappa(\gamma)$ is a protected guarantee for lower game with $c^U_{h,\gamma}$
  for all $h,\gamma>0.$

   Since $\varkappa(\gamma)$ tends to $0$ as $\gamma\to\gamma_*$,
  $U$ is a subsolution for lower game with payoff $\nu_\gamma(\gamma\to\gamma_*)$.
\bo
\begin{remark}\label{111}
 As follows from this proof,  under condition $(\tau)$,
  $\varkappa(\gamma)$-optimal strategies for payoff $\nu_{\gamma_h}$  protect the corresponding asymptotic guarantee for  games with $c^U_{h,\gamma}.$
\end{remark}

{\bf Proof of Lemma~\ref{firstlemm1}.}

 Let  $\mA$ satisfy conditions $(\ct{P}),(\omega)$. Let a mapping $c:\mm{K}\to \mm{R}$ be bounded.
 We must prove that, for all $\epsi>0$, $\cym[c]-\epsi$ is a protected guarantee for the lower game with payoff $c$.

 Indeed, for $\epsi>0$, for all $\omega\in\Omega$,
 there exists a $A^\omega\in\mA$ such that $c(z)>\cym[c](\omega)-\epsi$ for all $z\in A^\omega\cap\Gamma(\omega).$ Set $\eta:\Omega\to\mA$
 as follows: $\eta(\omega)=A^\omega$ for all $\omega\in\Omega.$ By condition $(\omega)$,
  there exists a strategy $A^*\in\mA$ such that $A^*\cap\Gamma(\omega)=\eta(\omega)\cap\Gamma(\omega)=A^\omega\cap\Gamma(\omega).$
  Then, $c(z)\geq\cym[c](\omega)-\epsi$ for all $z\in A^*\cap\Gamma(\omega),$ $\omega\in\Omega$, i.e.
  $\cym[c]-\epsi$ is a protected guarantee. \bo

{\bf Proof of Lemma~\ref{firstlemm2}.}

  Let  $\mA,\mB$ satisfy condition $(\ct{C})$, let the payoff $c:\mm{K}\to \mm{R}$ be bounded.
 We claim that $\cym[c]\leq\cyp[c]$.

Fix $\epsi>0$, $\omega\in\Omega.$  There exists $A^\omega\in\mA,B^\omega\in\Omega$ such that $c(z')>\cym[c](\omega)-\epsi$,
 $c(z'')<\cyp[c](\omega)+\epsi$ for all $z'\in A^\omega\cap\Gamma(\omega),z''\in B^\omega\cap\Gamma(\omega).$
 By condition $(\ct{C})$, there exists  $z\in A^\omega\cap B^\omega\cap\Gamma(\omega)$.
 Then,
 $\cym[c](\omega)-\epsi<c(z)<\cyp[c](\omega)+\epsi.$
 Since $\epsi>0,\omega\in\Omega$ were chosen arbitrarily, we get $\cym[c]\leq\cyp[c].$ \bo

\section{A guarantee for $\avm$ as a guarantee  for $\bwm$. The proof of Proposition \ref{av_bw}.}
\label{A}

 {\bf Step 1. Preliminary constructions and estimates.}

Notice that $U$
is bounded from above by a positive $R.$

It is easy to verify that $\ln p<p-1< p\ln p$ if $p>1.$ To each natural number $k>2$, we can assign a number
 $p\in(1,2)$  such that
 \begin{eqnarray}
\label{neqp}
\frac{1}{k}<\frac{\ln k}{k}
<\ln p<p-1< p\ln p<\frac{2 \ln k}{k}.
\end{eqnarray}
Fix such $k,p.$

Now, $S$ is a subsolution, therefore  for some $\hat{T}^{(k)}>0$, we have \rref{sol207}:
for all positive $T$,$\delta$ with $T>2\hat{T}^{(k)},\delta<T/2$, there exists  a strategy $A^{T,\delta}\in\mA$ such that,
for all $z\in A^{T,\delta}$,
\begin{eqnarray}
 U_{T}(z(0))-\frac{1}{k^2}&\leq&
  \frac{1}{T}\int_{0}^{\delta}g(z(t))\,dt+\frac{T-\delta}{T}U_{T-\delta}(z(\delta))\nonumber\\
  &\ravref{2000b}&
  \av_T(z)-\frac{T-\delta}{T}\Big(\av_{T-\delta}(z_{\delta})-U_{T-\delta}\big(z(\delta)\big)\Big).
 \label{sol207_}
\end{eqnarray}
   Moreover,  by \rref{slowlyT}, we also can choose $\hat{T}^{(k)}$ such that
   \begin{eqnarray}
\label{neqq}
 U_{T}(\omega)\geq U_{p^{-1}T}(\omega)-\frac{1}{k^2}\qquad \forall \omega\in\Omega, T>2\hat{T}^{(k)}.
\end{eqnarray}
  Fix such $T$.
Define
 \begin{eqnarray*}
 \lambda\rav\frac{1}{T},\quad
 \delta\rav\frac{T(p-1)}{p},\quad
  \tau_i\rav i\delta\qquad
  \forall i\in\overline{0,k}.
 \end{eqnarray*}
Then,
 \begin{eqnarray}
  \frac{T-\delta}{T}=p^{-1},\quad \frac{\delta}{T\ln p}=\frac{p-1}{p\ln p}\leqref{neqp} 1,
  \label{def1dob}
\\
\frac{p\delta}{T\ln p}\leqref{def1dob} p\leqref{neqp}
1+\frac{2\ln k}{k}=\frac{1}{\lambda T}+\frac{2\ln k}{k},&\ &\frac{\ln p}{\delta}\geqref{def1dob} \frac{1}{T}=\lambda.
\label{def1dob2}
 \end{eqnarray}

{\bf Step 2. Constructing a near-$\bw_\lambda$ payoff.}

 Define a piecewise constant function $\varrho$ on $[0,\tau_k)$ by the rule
\begin{eqnarray}
\label{neq16}
 \varrho(t)=p^{-i}\qquad \forall t\in [\tau_{i},\tau_{i+1}).
\end{eqnarray}
Then, for $t\in [0,\tau_k)$, we have $t\in[\tau_{i},\tau_{i+1})$ for some $i$, and
\begin{eqnarray}
\label{neq15}
\varrho(t)=p^{-i}= p^{1-\tau_{i+1}/\delta}\leq p^{1-t/\delta}.
\end{eqnarray}

Consider a lower  game with the following payoff:
\begin{eqnarray*}
{c}(z)\rav
\frac{1}{T} \int_{0}^{\tau_k} \varrho(t)g(z(t))dt+p^{-k} U_{p^{-1}T}(z(\tau_k))\qquad \forall z\in\mm{K}.
\end{eqnarray*}
 Note that, by  $U_{p^{-1}T}\leq R$,
 we have 
 $$p^{-k}U_{p^{-1}T}(\omega)\leq p^{-k} R=e^{-k\ln p} R\leqref{neqp}e^{-\ln k} R=\frac{R}{k}\leqref{neqp}\frac{R\ln k}{k}.
 $$
Now, by $0\leq g\leq 1$, for every process $z\in\mm{K}$,
\begin{eqnarray}
{c}(z)-\bw_\lambda(z)\!\!\!\!&\leqref{neq15}&\!\!\!\!
\frac{1}{T}\int_{0}^{\infty}\!\!\!\![p^{1-t/\delta}-e^{-\lambda t}]dt+\frac{R\ln k}{k}
\!\!=\!\!
\frac{ p\delta}{T\ln p}-\frac{1}{\lambda T}+\frac{R\ln k}{k}\leqref{def1dob2}\frac{(R+2)\ln k}{k}.
\label{neq56}
\end{eqnarray}

{\bf Step 3. Constructing the strategy $A^*$.}

 Remember that $T>2\hat{T}^{(k)}$, $\delta=(1-p^{-1})T<T/2$
 by choices of $T$ and $p$ respectively. Then, \rref{sol207_} holds for some $A=A^{T,\delta}\in\mA$.
 Let us also note that, since the right-hand side of this inequality depends only on $z|_{[0,\delta]}$,
 the strategy $A$  from \rref{sol207_} can be replaced with arbitrary strategy that could be represented in the form
  $A^T\diamond_{\delta}A.$
  Now, by \rref{def1dob}, it is equivalent to
  \begin{equation}
\label{1050}
   U_T(z(0))-\frac{1}{k^2}\leq
  \av_{T}(z)-p^{-1}
  \av_{p^{-1}T}(z_{\delta})+p^{-1} U_{p^{-1}T}(z(\delta))
   \qquad
   \forall A'\in\mA,z\in A\diamond_{\delta}A'.
 \end{equation}
 Define
 \begin{equation}
\label{UT}
 A^*=A\diamond_{\tau_1}A\diamond_{\tau_2}\dots\diamond_{\tau_{i-1}}A\diamond_{\tau_i}
 \dots\diamond_{\tau_{k-1}}A\diamond_{\tau_{k}}A.
\end{equation}
Note that such $A^*$ exists in $\mA$ by property $(\diamond)$.

  Note that, since the sufficiently large $k$ were chosen arbitrarily,
  to prove the proposition, it would suffice to demonstrate that
\begin{eqnarray*}
\bwm_{\lambda}(\omega)\geq U_{T}(\omega)-\frac{(R+4)\ln k}{k}\qquad \forall\omega\in\Omega;
\end{eqnarray*}
  in accordance with \rref{neq56}, this fact would follow from
\begin{equation}
\label{neq66}
{c}(z)>U_T(z(0))-\frac{2\ln k}{k}\qquad\forall z\in A^*.
\end{equation}

{\bf Step 4. Proof of estimate \rref{neq66}.}

 Remember that $p^{-1}T=T-\delta,$
                $\tau_{i+1}=\tau_{i}+\delta,$   
                $\varrho(t)=p^{-i}$ for $t\in[\tau_{i},\tau_{i+1}).$
                Now,
\begin{eqnarray*}
  \frac{1}{T}\int_{\tau_{i}}^{\tau_{i+1}} \varrho(t)g(z(t))dt&=&
   \frac{1}{T}\int_{0}^{\delta} \varrho(\tau_{i})g(z(t+\tau_{i}))dt\\&=&
   \frac{p^{-i}}{T}\int_{0}^{\delta} g(z(t+\tau_{i}))dt
   \ravref{2000b}
  p^{-i} \av_T(z_{\tau_{i}})-p^{-i-1} \av_{p^{-1}T}(z_{\tau_{i+1}}).
\end{eqnarray*}
  Then, for a process $z\in\mm{K},$ we obtain
\begin{eqnarray}
  {c}(z)= \av_{T}(z)-p^{-1}\av_{p^{-1}T}(z_{\tau_{1}})&+&
\dots\nonumber\\
  p^{-i}\av_{T}(z_{\tau_{i}})-
  p^{-i-1}\av_{p^{-1}T}(z_{\tau_{i+1}})&+&\dots\nonumber\\
  p^{-k+1} \av_{T}(z_{\tau_{k-1}})-
  p^{-k}
  \av_{p^{-1}T}(z_{\tau_{k}})&+&
  p^{-k}U_{p^{-1}T}(z(\tau_{k})).\label{last25}
\end{eqnarray}
 By \rref{UT}, for every $z\in A^*$, we have $z_{\tau_{k-1}}\in A\diamond_\delta A$; then,
 thanks to \rref{1050}, if we take into account $\tau_{k-1}+\delta=\tau_{k}$, we will obtain
\begin{eqnarray*}
  \av_{T}(z_{\tau_{k-1}})- p^{-1}
  \av_{p^{-1}T}(z_{\tau_{k}})+p^{-1} U_{p^{-1}T}(z(\tau_{k}))&\geqref{1050}&
  U_T(z(\tau_{k-1}))-\frac{1}{k^2}
  \geqref{neqq}
   U_{p^{-1}T}(z(\tau_{k-1}))-\frac{2}{k^2}.
\end{eqnarray*}
 Substituting this into \rref{last25} and accounting for $\tau_{k-1}+\delta=\tau_{k}$, we get
\begin{eqnarray}
 {c}(z)\geq \av_{T}(z)-p^{-1}\av_{p^{-1}T}(z_{\tau_{1}})&+&
  \dots\nonumber\\
  p^{-i}\av_{T}(z_{\tau_{i}})-
  p^{-i-1}\av_{p^{-1}T}(z_{\tau_{i+1}})&+&\dots\nonumber\\
  p^{-k+2} \av_{T}(z_{\tau_{k-2}})-
  p^{-k+1}  \av_{p^{-1}T}(z_{\tau_{k-1}})&+&p^{-k+1}U_{p^{-1}T}(z(\tau_{k-1}))-\frac{2}{k^2}
  \quad \forall z\in A^*.\label{last1_1}
\end{eqnarray}
 From \rref{UT} and $\tau_{k-2}+\delta=\tau_{k-1}$,
 we obtain $z_{\tau_{k-2}}\in A\diamond_{\delta}(A\diamond_\delta A)$
 for every $z\in A^*$. Now,
\begin{eqnarray*}
  \av_{T}(z_{\tau_{k-2}})- p^{-1}
  \av_{p^{-1}T}(z_{\tau_{k-1}})+p^{-1} U_{p^{-1}T}(z(\tau_{k-1}))&\geqref{1050}&
  U_T(z(\tau_{k-2}))-\frac{1}{k^2}
  \\ &\geqref{neqq}&
    U_{p^{-1}T}(z(\tau_{k-2}))-\frac{2}{k^2}.
\end{eqnarray*}
 Substituting this into \rref{last1_1}, we obtain
\begin{eqnarray*}
  {c}(z)\geq \av_{T}(z)-p^{-1}\av_{p^{-1}T}(z_{\tau_{1}})&+&\dots\\
  p^{-i}\av_{T}(z_{\tau_{i}})-
  p^{-i-1}\av_{p^{-1}T}(z_{\tau_{i+1}})&+&\dots\\
  p^{-k+3} \av_{T}(z_{\tau_{k-3}})-
  p^{-k+2}  \av_{p^{-1}T}(z_{\tau_{k-2}})&+&p^{-k-2}U_{p^{-1}T}(z(\tau_{k-2}))- \frac{4}{k^2}
    \quad \forall z\in A^*.
\end{eqnarray*}
 Proceeding in the similar way, since it is always $\tau_{k-l}+\delta=\tau_{k-l+1}$
 and
\begin{eqnarray*}
  \av_{T}(z_{\tau_{k-l}})- p^{-1}
  \av_{p^{-1}T}(z_{\tau_{k-l+1}})+p^{-1} U_{p^{-1}T}(z(\tau_{k-l+1}))\geq U_{p^{-1}T}(z(\tau_{k-l+1}))-\frac{2}{k^2}
\end{eqnarray*}
 holds for all $z\in A^*$, we now see
 that, for every instance of $z\in A^*$, it holds that
\begin{eqnarray*}
  {c}(z)\geq \av_{T}(z)-p^{-1}\av_{p^{-1}T}(z_{\tau_{1}})+p^{-1}U_{p^{-1}T}(z(\tau_{1}))- \frac{2k-2}{k^2}.
\end{eqnarray*}
 By $\tau_1=\delta$ and $\ln k>1,$ the relation \rref{1050} directly implies \rref{neq66} for all $z\in A^*$,
  which was to be proved. \bo
\begin{remark}\label{222}
 As follows from the proof, the strategies that protect  an
  asymptotic guarantee $U$ for  the payoffs $\bw_\lambda,$
  can be constructed by rule \rref{UT} with  the aid of the strategies $A^{T,h}$ that protect
  the similar asymptotic guarantee  for  payoffs  $\hat{c}^U_{T,h}.$
\end{remark}




\section{ A guarantee for $\bwm$ is a guarantee for $\avm$.  The proof of Proposition \ref{bw_av}}
\label{B}


 {\bf Step 1. Preliminary constructions and estimates.}

Notice that $U$
is bounded from above by a positive $R.$

Consider a natural $k$; we can map to it the numbers $M>1,p>1$  such that
 $$k=M\ln M,\qquad p\rav e^{1/M}.$$
 Note that $1+x<e^x<1+x+x^2$  for $|x|\in(0,1).$
 By $M>1$, we have $p=e^{1/M}=1+\frac{1}{M}+\frac{r'}{M^2}$,
 $p^{-1}=e^{-1/M}=1-\frac{1}{M}+\frac{r''}{M^2}$ for some $r',r''\in(0,1)$.
Then,
\begin{eqnarray}
\frac{1-\frac{p}{M}}{M(1-p^{-1})}=\frac{1-\frac{1}{M}-\frac{1}{M^2}-\frac{r'}{M^3}}{1-\frac{r''}{M}}<1.
\label{neqM}
\end{eqnarray}

  By  \rref{slowlyl},
for  $k$ (and, therefore, $M$)
there exists $T_{0}>k=M\ln M$ such that
\begin{eqnarray}
\label{neq1}
U_{\frac{pM}{T}}\leq U_{\frac{M}{T}}+\frac{1}{k^2}\qquad \forall T>T_0.
\end{eqnarray}
Since $U$ is a subsolution, thanks to \rref{sol307}, we also can choose $T_0$ such that,
  for each positive $\lambda<{M}/{T}_0,$ for each $h>0$,
 there exists $A^{\lambda,h}\in\mA$ such that, for all  $z\in A^{\lambda,h}$,
\begin{eqnarray}
 U_{\lambda}(z(0))-\frac{1}{k^2}\!<\!
  \lambda\int_{0}^{h}e^{-\lambda t}g(z(t))\,dt+e^{-\lambda h} U_{\lambda}(z(h))\!\!\ravref{2000b}\!\! \bw_\lambda(z)-e^{-\lambda h}\bw_\lambda(z_{h})+e^{-\lambda h} U_{\lambda}(z(h)).
 \label{sol307_}
\end{eqnarray}

 Let us fix such $k,M,p,T_0.$ Fix also some $T>T_{0}.$
Define
\begin{eqnarray*}
\lambda=\frac{1}{T},\ \quad
 t_0=\frac{T}{M},\ 
 \quad\tau_0=0,\ \quad
 t_i=t_0 p^{-i},\ \quad
 \tau_i=\tau_{i-1}+t_{i-1}&\ &\quad\forall i\in\overline{1,k}.
 \end{eqnarray*}
 Then, we have the following succession of equalities and inequalities:
\begin{eqnarray}
\label{neq3}
  p^{-k}=e^{-\ln M}=\frac{1}{M},
  &\ & \\
\label{neq4}\quad
\lambda p^i\leq\lambda p^k=\lambda M=\frac{M}{T}\leq\frac{M}{T_{0}},\quad
U_{\lambda p^i}(\omega)\geqref{neq1}U_{\lambda p^{i-1}}(\omega)-\frac{1}{k^2}&\ &\qquad
\forall i\in\overline{1,k},\omega\in\Omega.
\end{eqnarray}
 Let us also note that $t_i$ constitute a monotonically decreasing geometric progression;
 $\tau_i$ are their partial sums, and
\begin{eqnarray}
\label{neqtau}
 \frac{1-\frac{p}{M}}{M(1-p^{-1})}
 \ravref{neq3}\frac{1-p^{-k+1}}{M(1-p^{-1})}=\frac{\tau_k}{T}\leqref{neqM}1.
 \end{eqnarray}

{\bf Step 2. Constructing a near-$\av_T$ functional.}

 Define a scalar function~$\varrho$ on $(0,\tau_k]$ by
 $$\varrho(t)=e^{-\lambda p^i (t-\tau_{i-1})}\qquad \forall i\in \overline{1,k}, t\in (\tau_{i-1},\tau_{i}].$$
 Note that on each such subinterval,
 \begin{eqnarray}
\label{neq5}
1\geq \varrho(t)\geq
 e^{-\lambda p^{i} (\tau_{i+1}-\tau_{i})}=e^{-\lambda p^{i} t_i}=e^{-\lambda t_0}=e^{-1/M}=p^{-1}>
 1-\frac{1}{M},\quad\forall i\in \overline{1,k}, t\in (\tau_{i-1},\tau_{i}].
 \end{eqnarray}

Consider a lower game with the following payoff
\begin{eqnarray*}
{c}(z)\rav
\lambda \int_{0}^{\tau_k} \varrho(t)g(z(t))\,dt+p^{-k} U_{\lambda p^{k}}(z(\tau_{k}))=
\frac{1}{T} \int_{0}^{\tau_k} \varrho(t)g(z(t))\,dt+p^{-k} U_{\lambda p^{k}}(z(\tau_{k})),\quad \forall z\in\mm{K}.
\end{eqnarray*}

 Recall that $0\leq g\leq 1$, $U\leq R;$ now,
 for every process $z\in\mm{K},$
$$p^{-k}U_{\lambda p^{k}}(z(\tau_k))\leq p^{-k} R\ravref{neq3}\frac{R}{M}$$
implies
\begin{equation}
\label{neq55}
 \av_T(z)\geqref{neqtau}
 \frac{1}{T}\int_{0}^{\tau_k}g(z(t))\,dt \geqref{neq5}
 \frac{1}{T}\int_{0}^{\tau_k}\varrho(t)g(z(t))\,dt
\geq
 c(z)-\frac{R}{M} \qquad \forall z\in\mm{K}.
\end{equation}

{\bf Step 3. Constructing strategy $A^*$.}

Note that, for every 
$i=\overline{0,k-1}$, we have
$e^{-\lambda p^{i} t_i}=p^{-1}$ by \rref{neq5} and $\lambda p^{i}<M/T_0$ by \rref{neq4}. Then,
 by~\rref{sol307_},
 there exists a strategy $A^{(i)}\rav A^{\lambda p^i,t_i}\in\mA$ such that
\begin{eqnarray}
 U_{\lambda p^i}(z(0))\leq
   \bw_{\lambda p^i}(z)-p^{-1}
  \bw_{\lambda p^i}(z_{t_i})+p^{-1} U_{\lambda p^i}(z(t_i))+\frac{1}{k^2}\qquad\forall z\in A^{(i)}.
\end{eqnarray}
 Since the right-hand side of this inequality  depend only on $z|_{[0,t_i]}$,
 the strategy $A^{(i)}$  can be replaced with arbitrary strategy that could be expressed in the form
  $A^{(i)}\diamond_{t_i}A'.$
Thus,
  \begin{equation}
\label{1749}
   U_{\lambda p^i}(z(0))\leq
  \bw_{\lambda p^i}(z)-p^{-1}
  \bw_{\lambda p^i}(z_{t_i})+p^{-1} U_{\lambda p^i}(z(t_i))+\frac{1}{k^2}
   \quad \forall i\in\overline{0,k-1},A'\in\mA,z\in A^{(i)}\diamond_{t_i}A'.
\end{equation}

 By property $(\diamond)$, there exists a strategy
\begin{equation}
\label{Ul}
A^*\rav A^{(0)}\diamond_{\tau_1}A^{(1)}\diamond_{\tau_2}\dots
 \dots\diamond_{\tau_{k-1}}A^{(k)}\in\mA.
\end{equation}

In view of  $\tau_{i+1}=\tau_{i}+t_{i}$, for $A^*$  constructed in this way, there exists $A'\in\mA$ such that
$z_{\tau_i}\in A^{(i)}\diamond_{t_i} A'$
for all $z\in A^*,i\in\overline{0,k-1}$. Now,
\rref{1749} implies
$$   U_{\lambda p^i}(z(\tau_i))\leq
  \bw_{\lambda p^i}(z_{\tau_i})-p^{-1}
  \bw_{\lambda p^i}(z_{\tau_{i+1}})+p^{-1} U_{\lambda p^i}(z(\tau_{i+1}))+\frac{1}{k^2}
   \qquad \forall z\in A^{*},i\in\overline{0,k-1}.$$
Finally, from \rref{neq4}, we  obtain
  \begin{equation}
\label{1750}
  U_{\lambda p^i}(z(\tau_i))\leq
  \bw_{\lambda p^i}(z_{\tau_i})-p^{-1}
  \bw_{\lambda p^i}(z_{\tau_{i+1}})+p^{-1} U_{\lambda p^{i+1}}(z(\tau_{i+1}))+\frac{2}{k^2}
   \qquad \forall z\in A^{*},i\in\overline{0,k-1}.
 \end{equation}

 Recall that $A^*\in\mA$ and, for all $\omega\in\Omega,$ we have
 $$  \avm_T(\omega)\geq \inf_{z\in A^*\cap \Gamma(\omega)} \av_T(z).$$
   Since $M$ (with $k=M\ln M$) can be arbitrary large,
   to prove the proposition, it would suffice to prove the inequality
\begin{eqnarray*}
\av_T(z)>U_\lambda(z(0))-\frac{R}{M}-\frac{2}{M\ln M}
\qquad \forall z\in A^*,
\end{eqnarray*}
  which follows from \rref{neq55} and 
\begin{equation}
\label{neq6}
{c}(z)\geq U_\lambda(z(0))-\frac{2}{k}\qquad\forall z\in A^*.
\end{equation}

{\bf Step 4. Proof of estimate \rref{neq6}.}

 For each $z\in\mm{K}$,
   $i=\overline{0,k-1}$, one has
\begin{eqnarray*}
  \lambda\int_{\tau_{i}}^{\tau_{i+1}} \varrho(t)g(z_t)dt&=&
  \lambda\int_{\tau_{i}}^{\tau_{i+1}}e^{-\lambda p^{i} (t-\tau_{i})}g(z(t))\,dt\\
  &=&
  \lambda\int_{0}^{\infty} e^{-\lambda p^{i} t}g(z(t+\tau_i))\,dt
  -\lambda e^{-\lambda p^{i} (\tau_{i+1}-\tau_{i})}\int_{0}^{\infty} e^{-\lambda p^{i} t}g(z(t+\tau_{i+1}))\,dt\\
  &\ravref{neq5}&
  \lambda\int_{0}^{\infty} e^{-\lambda p^{i} t}g(z_{\tau_i}(t))\,dt
  -\lambda p^{-1}\int_{0}^{\infty} e^{-\lambda p^{i} t}g(z_{\tau_{i+1}}(t))\,dt\\
  &=&
  p^{-i} \bw_{\lambda p^{i}}(z_{\tau_{i}})-
  p^{-i-1} \bw_{\lambda p^{i}}(z_{\tau_{i+1}}).
\end{eqnarray*}
Then, for each $z\in A^*,$ $i=\overline{0,k-1}$, we have
\begin{eqnarray*}
  \lambda\int_{\tau_{i}}^{\tau_{i+1}} \varrho(t)g(z_t)\,dt&=&
  p^{-i} \bw_{\lambda p^{i}}(z_{\tau_{i}})-
  p^{-i-1} \bw_{\lambda p^{i}}(z_{\tau_{i+1}})+p^{-i-1} U_{\lambda p^{i+1}}(z(\tau_{i+1}))-p^{-i-1} U_{\lambda p^{i+1}}(z(\tau_{i+1}))\\
  &=&
  p^{-i}\Big( \bw_{\lambda p^{i}}(z_{\tau_{i}})-
  p^{-1} \bw_{\lambda p^{i}}(z_{\tau_{i+1}})+p^{-1} U_{\lambda p^{i+1}}(z(\tau_{i+1}))\Big)-p^{-i-1} U_{\lambda p^{i+1}}(z(\tau_{i+1}))\\
  &\geqref{1750}& p^{-i} U_{\lambda p^i}(z(\tau_i))-p^{-i-1} U_{\lambda p^{i+1}}(z(\tau_{i+1}))-\frac{2}{k^2}.
\end{eqnarray*}
  Summing over all the intervals $[\tau_i,\tau_{i+1}]$, we get
\begin{eqnarray*}
  {c}(z)&=&p^{-k} U_{\lambda p^{k}}(z(\tau_{k}))+\sum_{i=0}^{k-1}\lambda\int_{\tau_{i}}^{\tau_{i+1}} \varrho(t)g(z_t)\,dt\\
        &\geq&p^{-k} U_{\lambda p^{k}}(z(\tau_{k}))+\sum_{i=0}^{k-1}\Big[p^{-i} U_{\lambda p^i}(z(\tau_i))-p^{-i-1} U_{\lambda p^{i+1}}(z(\tau_{i+1}))-\frac{2}{k^2}\Big]\\
        &=&p^{-k} U_{\lambda p^{k}}(z(\tau_{k}))+U_{\lambda}(z(\tau_0))-p^{-k} U_{\lambda p^{k}}(z(\tau_{k}))-\frac{2}{k}\\
        &=&U_{\lambda}(z(0))-\frac{2}{k}\qquad\qquad\qquad\qquad\forall z\in A^*.
\end{eqnarray*}
 Thus, inequality  \rref{neq6} is proved.\bo
\begin{remark}\label{333}
 As follows from the proof,  the strategies
  protecting an asymptotic guarantee $U$ for  payoffs $\av_T$ can be
  constructed by rule \rref{Ul}  with   the aid of strategies $A^{\lambda,h}$ that protect
  the similar asymptotic guarantee for payoffs like  $\check{c}^U_{\lambda,h}.$
\end{remark}

 \section{Tauberian theorem for differential games}
\label{D}

 Consider a nonlinear system in $\mm{R}^n$ controlled by two players
\begin{equation}\label{sys}
 \dot{x}=f(x,a,b),\ x(0)\in \mm{R}^n,\ a(t)\in \mm{A},\ b(t)\in \mm{B}\qquad a.e.\ t\geq 0;
\end{equation}
here,~$\mm{A}$ and~$\mm{B}$ are compact metric spaces.

 Here and below, we assume  functions
 $f:\mm{R}^n \times \mm{A} \times \mm{B}\to\mm{R}^n$, $g:\mm{R}^n \times \mm{A} \times \mm{B}\to [0,1]$ are
\begin{enumerate}
  \item continuous;
  \item Lipschitz continuous in the state variable, namely, for a constant $L>0$,
  $$\big|\big|f(x,a,b)-f(y,a,b)\big|\big|+\big|g(x,a,b)-g(y,a,b)\big|\leq L\big|\big|x-y\big|\big|\qquad \forall x,y\in\mm{R}^n,a\in \mm{A},b\in \mm{B}.$$
\end{enumerate}
  Remember that $B(\mm{T},\mm{A})$ and $B(\mm{T},\mm{B})$ are the sets of all Borel measurable functions $\mm{T}\ni t\mapsto a(t)\in \mm{A}$ and $\mm{T}\ni t\mapsto b(t)\in \mm{B}$, respectively.
  Since the elements of both $\mm{A}$ and $\mm{B}$ are functions, not equivalence classes, hereinafter, in this section, all equivalences default to everywhere, not almost everywhere.

 Now, for each pair  $(a,b)\in B(\mm{T},\mm{A})\times B(\mm{T},\mm{B})$, for every initial condition  $x(0)=x_*\in\mm{R}^n$, system \rref{sys} generates the unique solution  $x(\cdot)=y(\cdot;x_*,a,b)$ defined for all $\mm{T}$.
 Denote by  $Y(x_*)$ the set of all such solutions with $x(0)=x_*$.

  Consider a set $\mm{X}\subset\mm{R}^n$ that is strongly invariant with respect to system \rref{sys}, i.e.,
   let $x(t)\in\mm{X}$ for all $t\in\mm{T},$ $x_*\in\mm{X}$, $x\in Y(x_*)$. Define
  $\mm{Y}\rav \cup_{x_*\in\mm{X}} Y(x_*).$

   Let us further assume the Isaacs' condition (also referred to as 'solvability of the small game'~\cite{ks}) holds, i.e.,
 \begin{subequations}
  \begin{equation}
 \label{sedlo}
  \max_{a\in \mm{A}}\min_{b\in \mm{B}}\big[  s\cdot f(x,a,b)+g(x,a,b) ] =
   \min_{b\in \mm{B}}\max_{a\in \mm{A}}\big[ s\cdot f(x,a,b)+g(x,a,b) \big]\quad \forall x,s\in\mm{R}^n.
\end{equation}
 Easy see that, for each positive function $\varrho:\mm{R}\to\mm{R}_{>0}$, it implies, for all $t\in\mm{R}$,
  \begin{equation}
 \label{sedlo_r}
  \max_{a\in \mm{A}}\min_{b\in \mm{B}}\big[s \cdot  f(x,a,b)+\varrho(t)g(x,a,b) ] =
   \min_{b\in \mm{B}}\max_{a\in \mm{A}}\big[s \cdot  f(x,a,b)+\varrho(t)g(x,a,b) \big]\quad \forall x,s\in\mm{R}^n.
\end{equation}
 \end{subequations}
     The goal of the first player is to maximize the payoff function while the task of the second is to minimize it.
    Our payoff functions are the following:  for each  positive $T,\lambda$, for all $(x,a,b)\in\mm{Y}\times\mm{A}\times\mm{B}$,
     \begin{eqnarray} \label{277}
 \av_T(x,a,b)\rav\frac{1}{T}\int_{0}^T g(x(t),a(t),b(t))\,dt,\qquad
 \bw_\lambda(x,a,b)\rav\lambda\int_{0}^\infty e^{-\lambda t} g(x(t),a(t),b(t))\,dt.
 \end{eqnarray}
        Note that \rref{sedlo_r} for $\varrho(t)\equiv \frac{1}{T}$ and
        $\varrho(t)\equiv \lambda e^{-\lambda t}$ becomes the Isaacs' condition  for the payoff functions $\av_T$ and $\bw_\lambda$, respectively.

  There are many ways to define a game and the sets of strategies for each player; for a very well made review encompassing a large number of formalizations, refer to~\cite[Subsect.14,15]{subb}.
 The Isaacs' condition not only provides the equality of lower and upper values;
      in addition, it makes value functions independent of
      formalization of strategies
      \cite{kss},\cite[Subsect.~14]{subb},\cite{Bardi}.
      For the definition of the value of the game, we can employ, for example the nonanticipating strategies (see \cite{RN},\cite{EK}).
       Let  $\ct{A},\ct{B}$ be the sets of all nonanticipating strategies for the first player 
       and second player respectively (see \cite{subb},\cite[Definition~VIII.1.1]{Bardi}).
      For all $\lambda,T>0$, define the value functions $\avg_T,\bwg_\lambda$ as follows: for  all $ x_*\in\mm{R}^n$
\begin{eqnarray*}
   \avg_T(x_*)&\rav&\sup_{Q\in\ct{A}}\inf_{b\in B(\mm{T},\mm{B})}
   \av_T(y(\cdot;x_*,Q(b),b),Q(b),b)=\inf_{Q\in\ct{B}}\sup_{a\in B(\mm{T},\mm{A})}
  \av_T(y(\cdot;x_*,a,Q(a)),a,Q(a)),\\
   \bwg_\lambda(x_*)&\rav&\sup_{Q\in\ct{A}}\inf_{b\in B(\mm{T},\mm{B})}
   \bw_\lambda(y(\cdot;x_*,Q(b),b),Q(b),b)
   =\inf_{Q\in\ct{B}}\sup_{a\in B(\mm{T},\mm{A})}
  \bw_\lambda(y(\cdot;x_*,a,Q(a)),a,Q(a)).
\end{eqnarray*}

A variant of the following theorem
 was also announced in \cite{Khlopin1,Khlopin2} and proved in \cite{Khlopin3}
  with the aid of nonanticipating operators \cite{chdu},\cite{chs}.
\begin{theorem}
     \label{maintheoremdiff}
 Let $\mm{X}$ be strongly invariant with respect to system \rref{sys}.
 Assume Isaacs' condition \rref{sedlo}.

  The following limits exist, are uniform in $x_*\in\mm{X}$, and coincide
   $$ \lim_{T\uparrow\infty}\avg_T(x_*)
       =\lim_{\lambda\downarrow 0}\bwg_\lambda(x_*)\quad\forall x_*\in\mm{X}$$
if at least one of these limits exists and is uniform   in
 $x_*\in\mm{X}.$
    \end{theorem}
The proof of this theorem will appear to be  an immediate consequence of Theorem~\ref{maintheoremgeneral}.
    Since we only require the closedness of the set of strategies with respect to concatenation, we will choose among
     the  feedback-type formalizations.
    To make the references more convenient, take  the class of feedback  strategies  with
perfect  memory  and  perfect  state  measurement
 \cite[Ch. XIV]{ks};
see also \cite[Definition~VIII.3.1]{Bardi},\cite{elliott},\cite[Sect.~11]{kss}.
 For proof of identity between the values defined through  nonanticipating strategies and through feedback strategies with perfect memory and perfect state measurement, refer to
 \cite[Subsect.~14]{subb},\cite[Theorem VIII.3.11]{Bardi}.


{\bf Definition of feedback  $MM$-strategies.}\

\begin{definition}
A map  $\zeta: \mm{Y}\to B(\mm{R}_{>0},\mm{A})$  is { \it  a feedback MM-strategy}  (feedback  strategy with
perfect  memory  and  perfect  state  measurement \cite[Definition~VIII.3.1]{Bardi})  for the  first  player if\\
1) for  each  $t  >   0$,   $x|_{[0,t]}   =  y|_{[0,t]}$   implies $\zeta[x]|_{(0,t]}  =  \zeta[y]|_{(0,t]}$;\\
2) for  all  $x_*\in \mm{X},b\in B(\mm{T},\mm{B}), T>0$,  there  exists  a  unique Carath\'{e}odory solution $x(\cdot)=y(\cdot,x_*,\zeta,b)$ of
\begin{equation}
\label{300}
\dot{x}(t)=f(x(t),\zeta[x(\cdot)](t),b(t)),\qquad x(0)=x_*,\qquad a.e.\ t\in(0,T].
\end{equation}
\end{definition}

We denote  by $\fr{F}$ the  set  of feedback
MM-strategies  for  the  first  player.
This definition directly implies that a solution  $y(\cdot,x_*,\zeta,b)$ (on $(0,T]$)  can be uniquely extended up to a solution of \rref{300} (for a.a. positive $t$) from $\mm{Y}\subset C(\mm{T},\mm{X})$; therefore, we may take $y(\cdot,x_*,\zeta,b)\in\mm{Y}$.
 Now, denote $\alpha[x_*,\zeta,b]\rav \zeta[y(\cdot,x_*,\zeta,b)]\in B(\mm{R},\mm{A})$
  for all $x_*\in\mm{X},\zeta\in\fr{F},b\in B(\mm{T},\mm{B})$. This definition is well-defined; moreover,
  $$y(\cdot,x_*,\zeta,b)=y(\cdot,x_*,\alpha[x_*,\zeta,b],b).$$

Feedback  MM-strategies   for the  second  player  are introduced in the similar way.
 Assign the  set $\fr{G}$ of feedback
MM-strategies  to  the second player; for each
$\xi\in\fr{G}$, for each $x_*\in\mm{X}$, for each $a\in B(\mm{T},\mm{A})$,
there exist a unique  $y(\cdot,x_*,a,\xi)\in\mm{Y},$
 a unique
$\beta[x_*,a,\xi]\in B(\mm{T},\mm{B})$.

Thanks to
 the Isaacs condition,
 by \cite{ks},\cite{subb}  for
 payoffs $\av_T$,
 by  \cite[Theorem VIII.3.11]{Bardi} for payoffs $\bw_\lambda,$
 the values of upper and lower games coincide; moreover, these values coincide with the  values defined by
 nonanticipating strategies. Thus,  for all $T>0,\lambda>0,x_*\in\mm{R}^n$, we obtain
\begin{eqnarray*}
   \avg_T(x_*)&=&\sup_{\zeta\in\fr{F}}\inf_{b\in B(\mm{T},\mm{B})}
   \av_T(y(\cdot;x_*,\zeta,b),\alpha[x_*,\zeta,b],b)
   =\inf_{\xi\in\fr{G}}\sup_{a\in B(\mm{T},\mm{A})}
   \av_T(y(\cdot;x_*,a,\xi),a,\beta[x_*,a,\xi]);\\
   \bwg_\lambda(x_*)&=&\sup_{\zeta\in\fr{F}}\inf_{b\in B(\mm{T},\mm{B})}
   \bw_\lambda(y(\cdot;x_*,\zeta,b),\alpha[x_*,\zeta,b],b)
   =\inf_{\xi\in\fr{G}}\sup_{a\in B(\mm{T},\mm{A})}
   \bw_\lambda(y(\cdot;x_*,a,\xi),a,\beta[x_*,a,\xi]).
\end{eqnarray*}
 In particular, like in Sect.~\ref{abstractgame}, we may now assume that one player announces his  own feedback MM-    strategy  (from either
    $\fr{F}$, or $\fr{G}$, respectively) and another, knowing it, selects a measurable control
    (either from $B(\mm{T},\mm{B})$ or $B(\mm{T},\mm{A})$, respectively).

{\bf Proof of Theorem~\ref{maintheoremdiff}. Reduction to the abstract formulation.}\

    Since the function $g$ depends on $a,b$ in addition to depending on $x$, let us set
    $$\Omega\rav\mm{X}\times\mm{A}\times\mm{B},\quad \mm{K}\rav\{(y,a,b)\,|\,y\in\mm{Y},a\in B(\mm{T},\mm{A}),b\in B(\mm{T},\mm{B})\}.$$
    Moreover,
    the mappings   $\av_T,\bw_\lambda$ defined in view of \rref{277} correspond to those defined in view of  \rref{248}.
     Still, $\Gamma(\omega)=\{z=(x,a,b)\in\mm{K}\,|\,z(0)=\omega\}$ for all $\omega\in\Omega.$    It remains to describe $\mA,\mB.$

For each $\zeta\in\fr{F}$,
it is valid to define
\begin{eqnarray*}
  A_\zeta\rav\bigcup_{a\in B(\mm{T},\mm{A}),b\in B(\mm{T},\mm{B}),x_*\in \mm{X}}\!\!\!\!\!\!\!\!\{(x,a,b)\in\mm{K}\,|\,\ x\rav y(\cdot;x_*,\zeta,b),a|_{(0,\infty)}=\zeta[x]\};\
  \mA\rav\{A_\zeta\subset\mm{K}\,|\,\zeta\in\fr{F}\}.
\end{eqnarray*}
We can define $B_\xi$ for each $\xi\in\fr{G}$ and the set $\mB$ in a similar way.

 For all $T>0,\lambda>0$, define the mappings $\avm_T,\avp_T,\bwm_\lambda,\bwp_\lambda$
 from $\Omega$ to $[0,1]$
  by formulas \rref{280a},\rref{280b}.     Note that, since
  $y(\cdot;x_*,\alpha(b),b)$ is independent  of $a(0)$ and $b(0)$, each of payoffs $\av_T$,$\bw_\lambda,$
    $\check{c}^U_{h,\lambda}$,$\hat{c}^U_{h,T}$ (for all $h,T,\lambda>0,$ $U:\Omega\to[0,1]$)
      is independent of them as well.

{\bf Proof of Theorem~\ref{maintheoremdiff}. Verification of  conditions of Theorem~\ref{maintheoremgeneral}.}

We claim that, for all $\omega=(x_*,a_*,b_*)\in\Omega$, $T,\lambda>0$
\begin{equation}
 \label{1066}
 \avm_T(x_*,a_*,b_*)=\avp_T(x_*,a_*,b_*)=\avg(x_*),\  \bwp_\lambda(x_*,a_*,b_*)=\bwm_\lambda(x_*,a_*,b_*)=\bwg_\lambda(x_*).
\end{equation}
Indeed, for each choice $\zeta\in\fr{F}$ by the first player,
a choice $b\in B(\mm{T},\mm{B})$ by the second player determines the unique  process
$(x,a',b')$ with $x\rav y(\cdot;x_*,\zeta,b)$, $b'(0)=b_*,b'|_{(0,\infty)}=b|_{(0,\infty)}$,$a'(0)=a_*,a'|_{(0,\infty)}=\zeta[x](b).$
By   $y(\cdot,x_*,\zeta,b)=y(\cdot,x_*,\alpha[x_*,\zeta,b],b),$
$\alpha[x_*,\zeta,b]=\zeta[y(\cdot,x_*,\zeta,b)],$
this is equivalent to the choice of $A_\zeta\in\mA$ followed by the choice of  $z\in A_\zeta\cap\Gamma(\omega).$
 Thus, \rref{1066} holds. Moreover, now, $A_\zeta\cap\Gamma(\omega)\neq\pust$. This  implies
 condition $(\ct{P})$ for $\mA$.
 One can,  totally analogously, prove  condition $(\ct{P})$ for $\mB$.

We claim that  condition $(\diamond)$ holds for the introduced families $\mA,\mB$.
In view of the symmetry, we will only prove this fact for  $\mA.$ Fix $\zeta',\zeta''\in\fr{F},\tau>0.$
Let us define the mapping $\zeta:\mm{Y}\to B(\mm{R}_{>0},\mm{A})$ piecemeal: first
$\zeta[x]|_{(0,\tau]},$ then, $\zeta[x]|_{(\tau,\infty)}$.

Assume $\zeta[x]|_{(0,\tau]}=\zeta'[x]|_{(0,\tau]}$  for all $x\in\mm{Y}.$ Conditions 1),2) from the definition of MM-strategy
hold for $\zeta$ under positive $T\leq\tau$ because they hold for $\zeta'.$

Similar to \rref{296}, for all $\tau>0,x\in\mm{Y},b\in B(\mm{T},\mm{B})$ , define  $x_\tau,b_\tau$ as follows:
$x_\tau(t)=x(t+\tau)$,$b_\tau(t)=b(t+\tau)$ for all $t\in\mm{T}.$
Note that, since $f, \mm{B}$ are independent of $t$, we have $x_\tau\in\mm{Y},b_\tau\in B(\mm{T},\mm{B}).$

Then, we can define $\zeta[x](t+\tau)=\zeta''[x_\tau](t)$ for all $t>\tau,x\in\mm{Y},$ i.e., $(\zeta[x])_\tau=\zeta''[x_\tau].$
It is easy see that this mapping $\mm{Y}\ni x\mapsto\zeta[x]|_{(\tau,\infty)}$ is nonanticipating because
$\zeta''$ is  nonanticipating.
Therefore, the map $\zeta$ is nonanticipating.
To prove Condition~2) from the definition of MM-strategy for all $T>\tau$, note that,
for an initial condition $x_*\in\mm{X}$ at $0$ and arbitrary  $b\in B(\mm{T},\mm{B})$, there exists a unique solution $y(\cdot,x_*,\zeta,b)|_{[0,\tau]}$ of \rref{300}, in particular, $y(\tau,x_*,\zeta,b)$  is also well-defined. In addition,
because $\zeta''$ is an MM-strategy,  for an initial position  at $\tau$ (in particular, for position $y(\tau,x_*,\zeta,b)$), for  a control of second player, for every $T>\tau$,  there  exists  a  unique  solution  of \rref{300} on the interval $[\tau,T].$ Then, 2) holds for~$\zeta.$ Thus, $\zeta$ is an MM-strategy.

Now,   we have
$\alpha[x_*,\zeta,b]|_{[0,\tau]}=\alpha[x_*,\zeta',b]|_{[0,\tau]}$,
$\alpha[x_*,\zeta,b](t+\tau)=\alpha[y(\tau,x_*,\zeta,b),\zeta'',b_\tau](t)$
for all $x_*\in\mm{X}$,   $b\in B(\mm{T},\mm{B}),$ $t>\tau;$
i.e.,  $\{(x,a,b)\in A_\zeta\,|\,x(0)=x_*\}=\{(x,a,b)\in A_{\zeta'}\diamond_\tau A_{\zeta''}\,|\,x(0)=x_*\}.$
So, $A_\zeta=A_{\zeta'}\diamond_\tau A_{\zeta''}\in\mA$ for all $\zeta',\zeta''\in\fr{F},\tau>0$, and condition $(\diamond)$   holds for $\mA$.

Consider a  payoff $c$ among $\av_T$,$\bw_\lambda,$
    $\check{c}^U_{h,\lambda}$,$\hat{c}^U_{h,T}$ (for all $h,T,\lambda>0,$ $U:\Omega\to[0,1]$).
  This payoff $c$ is bounded and independent  of $a(0)$ and $b(0)$, therefore,
  by definition of the lower values,
for each $\epsi>0,$
  for each initial $x_*$,  there exists $\zeta^{x_*}$ with $\displaystyle\inf_{z\in A_{\zeta^{x_*}}} c(z)>\cym[c](z(0))-\epsi.$
  Define $\zeta_*:\mm{Y}\to B(\mm{T},\mm{B})$ as follows: $\zeta_*[y]\rav \zeta^{y(0)}[y]$ for all $y\in\mm{Y}.$
  By straightforward verification of  conditions 1)-2), it is easily proved that  $\zeta_*$ is MM-strategy.
  By the construction, this MM-strategy is an $\epsi$-optimal  MM-strategy for the payoff $c$.
  Therefore, for any $\epsi>0$, the lower game with this payoff  has $\epsi$-optimal strategy.
  The proof of existence of $\epsi$-optimal strategy for the upper game with this payoff
  for all $\epsi>0$ is analogous. It remains to prove that for all games mentioned in the formulation of  Theorem~\ref{maintheoremgeneral} there exists a saddle point.

Remember that,  by \cite[Subsect.~14]{subb}),\cite[Theorem VIII.3.11]{Bardi},
 for all payoffs
$\av_T,\bw_\lambda$ $(\lambda,T>0)$, the condition \rref{sedlo_r} for corresponding $\varrho$ implies
 that a differential game with dynamics \rref{sys} and this payoff has a saddle point.
 In particular, its value functions are bounded (lies in $[0,1]$) and   continuous \cite[Proposition VIII.1.8]{Bardi},\cite[Theorem~11.4]{subb}.
 By  \cite[Theorem~11.4]{subb}, it is the same for the payoff $$\int_{0}^h \varrho(t)g(z(t))\,dt+S(z(h))$$
 if the functions $S:\mm{R}^n\to\mm{R}$ and $\varrho:[0,h]\to\mm{R}$ are bounded and continuous.
  Then, it is the same for payoffs $\check{c}^S_{h,\lambda},\hat{c}^S_{h,T}$ for all $h,\lambda,T>0$
  and for each $S=\avm_T,S=\avp_T,S=\bwm_\lambda,S=\bwp_\lambda.$
So, each of the games needed for Theorem~\ref{maintheoremgeneral}
has a saddle point.

 All conditions of Theorem~\ref{maintheoremgeneral} hold. Applying this theorem, we prove
 the result of Theorem~\ref{maintheoremdiff}.
\bo
\end{document}